\documentclass{article}
\usepackage[utf8]{inputenc}
\usepackage{pgf,tikz,amsmath, amsfonts}

\begin{document}
\title{Mathematical specification of hitomezashi  designs}

\author{Katherine A.~Seaton$^{a}$\thanks{Email: k.seaton@latrobe.edu.au
} \ and Carol Hayes$^{b}$\\ $^{a}$Department of Mathematical and Physical Sciences,\\ La Trobe University VIC 3086, Australia\\
$^{b}$School of Culture, History and Language,\\ The Australian National University ACT 2600, Australia
 }

\maketitle

\begin{abstract} 
  Two mathematical aspects of the centuries-old Japanese sashiko stitching form hitomezashi are discussed: the encoding of designs using words from a binary alphabet, and duality. Traditional hitomezashi designs are analysed using these two ideas. Self-dual hitomezashi designs related to Fibonacci snowflakes, which we term Pell persimmon polyomino patterns, are proposed. Both these designs and the binary words used to generate them appear to be new to their respective literatures.


\textbf{AMS classes:}   97M80; 05A05

\end{abstract}

\section{Introduction}

The Japanese word {sashiko} literally means `little stabs', referring to the act of pushing one's needle vertically through thick layers of cloth. In the contemporary world of global stitching practice there are basically four types of sashiko stitching in use: {kogin/hishi/sh\={o}nai sashiko}  from the northern T\={o}hoku region of Japan,  {borozashi},  {hitomezashi} and {moy\={o}zashi}. All four stitching types are grounded in the cultural traditions of Japan dating back to the early Edo period, when sashiko was used to quilt together layers of roughly woven fabric, to mend well-used garments by adding strength and longevity and, finally, to decorate or embellish clothing and household textiles to turn the functional into the decorative. Hayes (2019) discusses the practical and cultural imperatives that led to the development of sashiko, how with societal change over time it fell out of everyday use, and then its revival as artform, nostalgic reference and mindfulness practice.  One contemporary purpose in the practice of sashiko is to ornament, decorating clothing and homewares, generally working on a single layer of fabric (Briscoe, 2004; \textit{Simply Sashiko}, 2019/2020). A second is practical, mending garments as part of the sustainable slow fashion movement (Marquez, 2019; Iiduka 2019/2021). 

Historically, the indigo-dyed cloth and cotton thread available to the poorer regions of northern Japan gave a distinctive blue-and-white look to sashiko garments and textiles. A notable feature of sashiko, compared, say, to needlepoint, is the deliberate use of empty space in designs, reflecting a distinctively Japanese aesthetic (Mende \& Morishige, 1991).

Tsugaru koginzashi and Nanbu {hishizashi}  both developed in the Aomori prefecture within T\={o}hoku, while sh\={o}nai sashiko  comes from the Yamagata prefecture, also in the north. All three of these types can be referred to as `diamond' embroidery in English as they use counted-thread stitching either on the horizontal or vertical to create a huge variety of rhombic patterns, with the stitches ranging in length from one to seven fabric threads (Hachinohe Institute of Technology, 2022). While kogin uses white thread, hishizashi commonly uses coloured thread, both in cotton and wool. The resulting geometric patterns tend to be drawn from activities of everyday life or images from the natural world, such as dragonflies or the diamond-leafed {hishi} (water chestnut). Sh\={o}nai sashiko is finer and more decorative, used by stitchers to express their wishes and prayers for a good harvest, for protection against evil, and for prosperous business. Flower patterns ({hanazashi}) carry thoughts of spring, and {sugizashi} (cedar stitch) evokes a sense of cedar covered mountains against the snow-covered fields of the northern winters. Contemporary stitchers using this style replicate the old patterns of the T\={o}hoku region.

Borozashi refers to the practice of using sashiko running stitch to mend and quilt together older and worn garments. The word boro translates as `scraps' or `shabby tattered rags' (Hayes, 2019). These scraps of material were layered onto older clothing, bedding or homewares and then quilted back into the article to extend its life. 

Moy\={o}zashi or pattern sashiko, now popular among stitchers globally and what many people think of when hearing the word sashiko,  refers to a contemporary sashiko stitching that is not restricted to traditional geometric patterning and stitch-length requirements. This style uses a series of dashed running stitches that never touch but can go in any direction, or trace curved lines, rather than following the traditional horizontal and vertical or diamond patterning. 
 In moy\={o}zashi the whole design is pre-printed or transferred to the fabric using a template (see, for example, Briscoe (2004)). 

Hitomezashi literally means `one-stitch sashiko', and refers to sashiko that is completed by stitching regularly on a grid, one running stitch per grid spacing. In a previous workshop paper we gave a (page-limited) outline of some of the mathematics that hitomezashi displays (Hayes \& Seaton, 2020). It is our intention in this article to discuss two particular mathematical aspects in detail: binary encoding of regular hitomezashi patterns, both traditional and what appear to be new patterns designed using mathematics, and duality.

It is the `one-stitch' nature of hitomezashi that permits the simple binary encoding of patterns outlined and used in this paper. Whereas traditionally the visible warp and weft of rough hemp fabric provided a square grid for these stitches, contemporary stitchers can also use specialist fabrics, such as even-weave linen or Aida cloth (see Figure 1). Since it is customary to stitch all the vertical lines of stitches and then all the horizontal lines, from edge to edge (unlike, for example, cross-stitch which is generally worked from the centre out), the design only becomes apparent as the second lot of stitches go in (see Figure \ref{kuchi}). This, too, encourages mindfulness. Hitomezashi can incorporate diagonal lines of stitches, or threads woven through the worked stitches; either of these is worked last.

\begin{figure}
\begin{center}
{{\includegraphics[height=6cm]{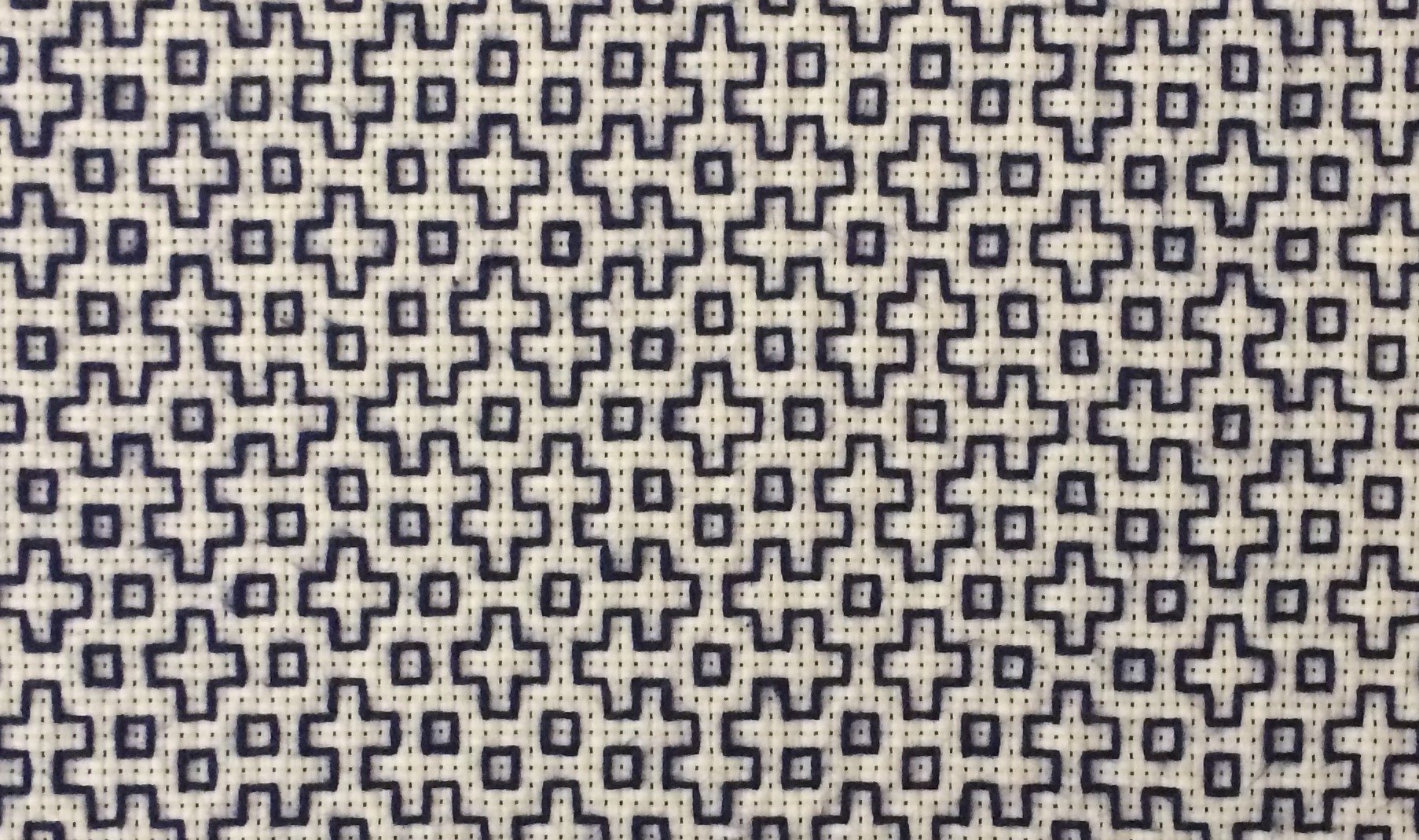}}}
\caption{\label{part} A mathematically-designed stitch pattern for hitomezashi. This piece has been worked on  white Aida cloth using blue cotton thread. The traditional roles of blue and white --- indigo-dyed cloth and white thread --- have been exchanged.}
\end{center}
\end{figure}

\begin{figure}
\begin{center}
{{\includegraphics[height=5cm]{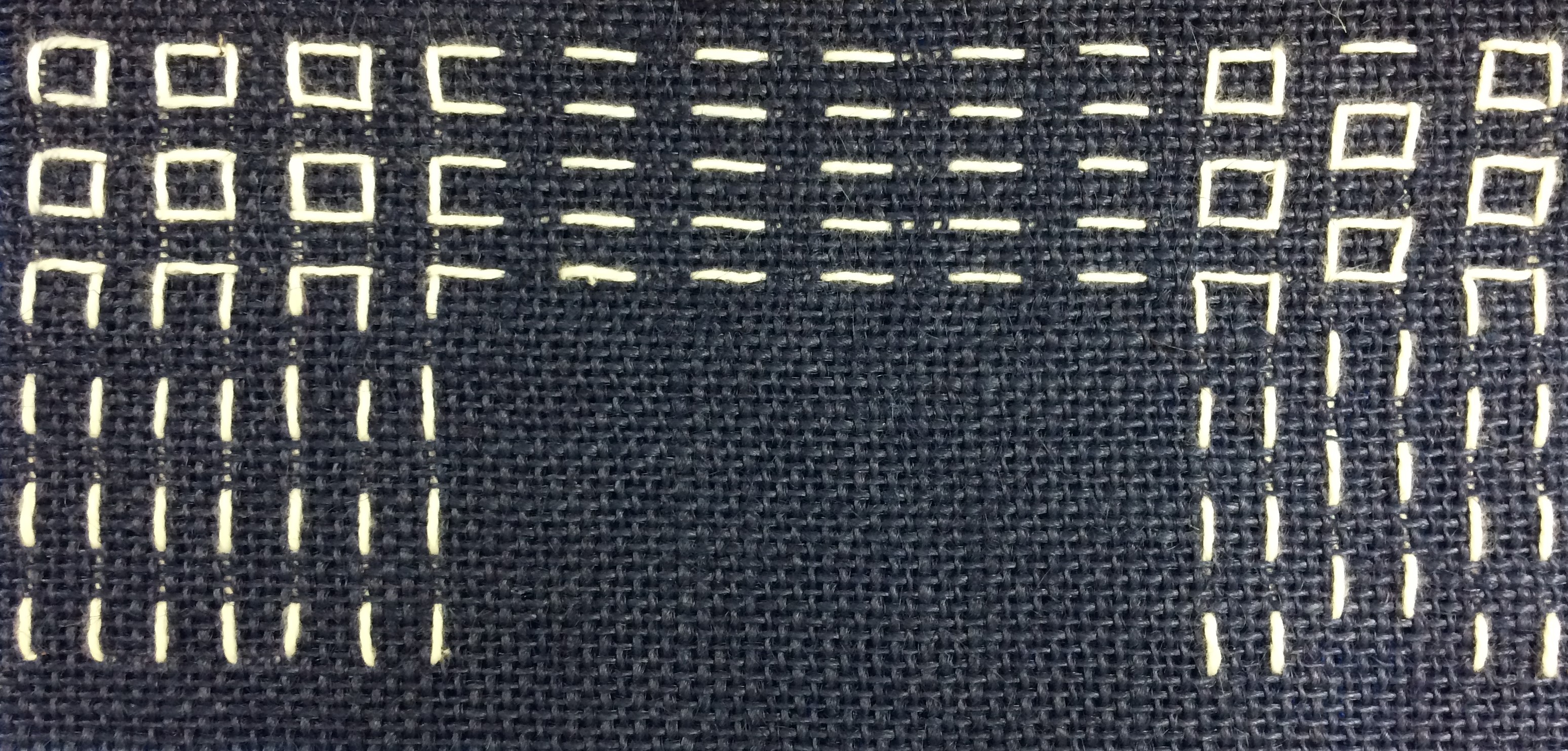}}}
\caption{\label{kuchi} Notice how different sequences of vertical stitches on the left and on the right interact with common lines of horizontal stitches to give different patterns. This piece has been worked with cotton thread on blue hessian (burlap). }
\end{center}
\end{figure}

A method of specifying hitomezashi designs in shorthand form using binary strings was employed by the first author when she proposed an activity for the Math Art Challenge organised by Annie Perkins (Perkins, 2020). The premise of the Challenge was to provide, each day for 100 days, an activity which could be done during lockdown or remote learning with materials on hand.  Hence it was suggested that hitomezashi designs could be drawn on graph paper, if sewing materials were not available. (Some participants chose to render them using their preferred software.) Three ways to choose the binary strings were suggested: to intentionally create regular patterns by repetition (which we explore in this paper), by flipping a coin (thus making aleotoric art) or as a form of steganography (e.g., by using the ascii representation of letters). 

The Math Art Challenge hitomezashi activity created a lot of interest on Twitter; subsequently, a Numberphile video featuring Ayliean MacDonald took the drawing idea to an even wider audience and explored the steganographic and aleotoric aspects (Haran, 2021). This video came to the attention of Defant and Kravitz (2022), who have subsequently proved a number of results about arbitrary hitomezashi designs, some of which we will refer to in Section \ref{loops}.

One observation made by many participants in the Math Art Challenge was that hitomezashi designs can be two-coloured. That is, only two colours are required to colour the regions marked out by the designs in such a way that adjacent regions have a different colour. Defant and Kravitz (2022) give the mathematical argument for this. The key observation is that hitomezashi stitching every vertical and horizontal line creates a design that is fully packed, that is, one in which every vertex has degree two.
While we had hinted at connections to statistical mechanics (Hayes \& Seaton, 2020), Defant and Kravitz (2022) identified specifically that results of Pete (2008) obtained in the context of corner percolation are actually results about square grid hitomezashi.  

More recently, Defant et al. (2022) have considered hitomezashi-like designs using stitches of different length in the vertical and horizontal directions and on the front and back of the work, on the square lattice. One observation that they make of relevance to this paper is that `normal' hitomezashi permits the greatest variety of patterns. They also consider the triangular lattice.

\begin{figure}
\begin{center}
{{\includegraphics[height=3cm]{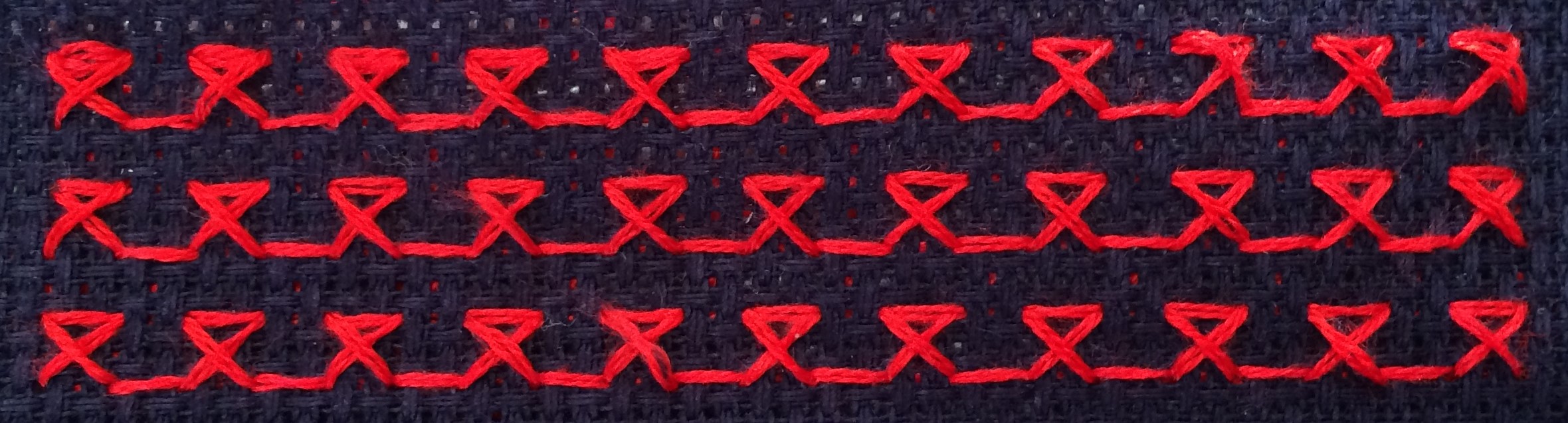}}}
\caption{\label{ribbons} Hitomezashi as a stitch form can incorporate crossing stitches and diagonal stitches. This design features both.  }
\end{center}
\end{figure}

In this article, as in our previous work (Hayes \& Seaton, 2020; Seaton, 2021) we confine our attention to hitomezashi worked in thread of one colour with horizontal and vertical lines of running stitch aligned to a square grid, such as that in Figure \ref{part}. With only a few exceptions, we consider designs that are fully packed.  We do not consider designs such as those of Figure \ref{ribbons} where there are diagonal lines or crossing stitches, or the {kugurizashi} form of hitomezashi in which threaded decoration is added. However, it is tiresome to repeat the restriction over and over again, so we state it here only.  We do want to acknowledge that the designs we consider comprise a subset of all the needlework patterns termed hitomezashi.

In the next section of the paper, we explain how hitomezashi patterns can be encoded using two symbols, and how duality is manifested when they are stitched (not drawn). We then discuss traditional hitomezashi stitch patterns explaining their evocative names, giving an encoding for each, and exploring their duality properties. In the third section of the paper we demonstrate that closed loops occurring in traditional patterns correspond to members of a family of polyominoes known as Fibonacci snowflakes. Equipped with this understanding of duality and the encoding, we propose new designs for hitomezashi (new in the sense that they do not appear in stitch dictionaries of traditional designs). We conjecture the relationship that these designs, which we term Pell persimmon polyomino patterns, or Pell persimmons for short, have to the Fibonacci snowflakes. We conclude by suggesting future work for ourselves and for our readers.

\section{Hitomezashi stitch patterns}

\subsection{Encoding the patterns}\label{code}
Making garments using sashiko to soften, strengthen or mend fabric was an essential domestic activity, particularly conducted in the winter when outside work was not possible (Hayes, 2019). Patterns were learnt from elders, by example and imitation (Briscoe, 2004). In contemporary sashiko books either written for non-Japanese speakers such as that of Briscoe (2004) or Marquez (2019), or translated for them such as \textit{Simply Sashiko} (2019/2020), a hitomezashi stitch design is generally presented in a chart showing a basic unit to be repeated to make a larger piece. 

Knitting patterns frequently use charts, but also give coded written instructions, in particular by using  K and P for the two basic stitches, knit and purl. (An excellent account of how knitting instructions have developed from descriptions given in words, how they use charts, and how they resemble coding languages is given by Howard (2017).) To give a simple example which motivates the encoding we will give for hitomezashi designs,  knitting patterns will frequently commence like this:
\smallskip

\noindent Cast on a multiple of six stitches.\\
Row 1 *K4\ P2; repeat from * to the end of the row.
\smallskip

\noindent We could say mathematically that the instruction for the first row is encoded by concatenation of the word KKKKPP from the binary alphabet $\{\text{K}, \text{P}\}$.

\begin{figure}
\begin{center}
{\resizebox*{6cm}{!}{\includegraphics{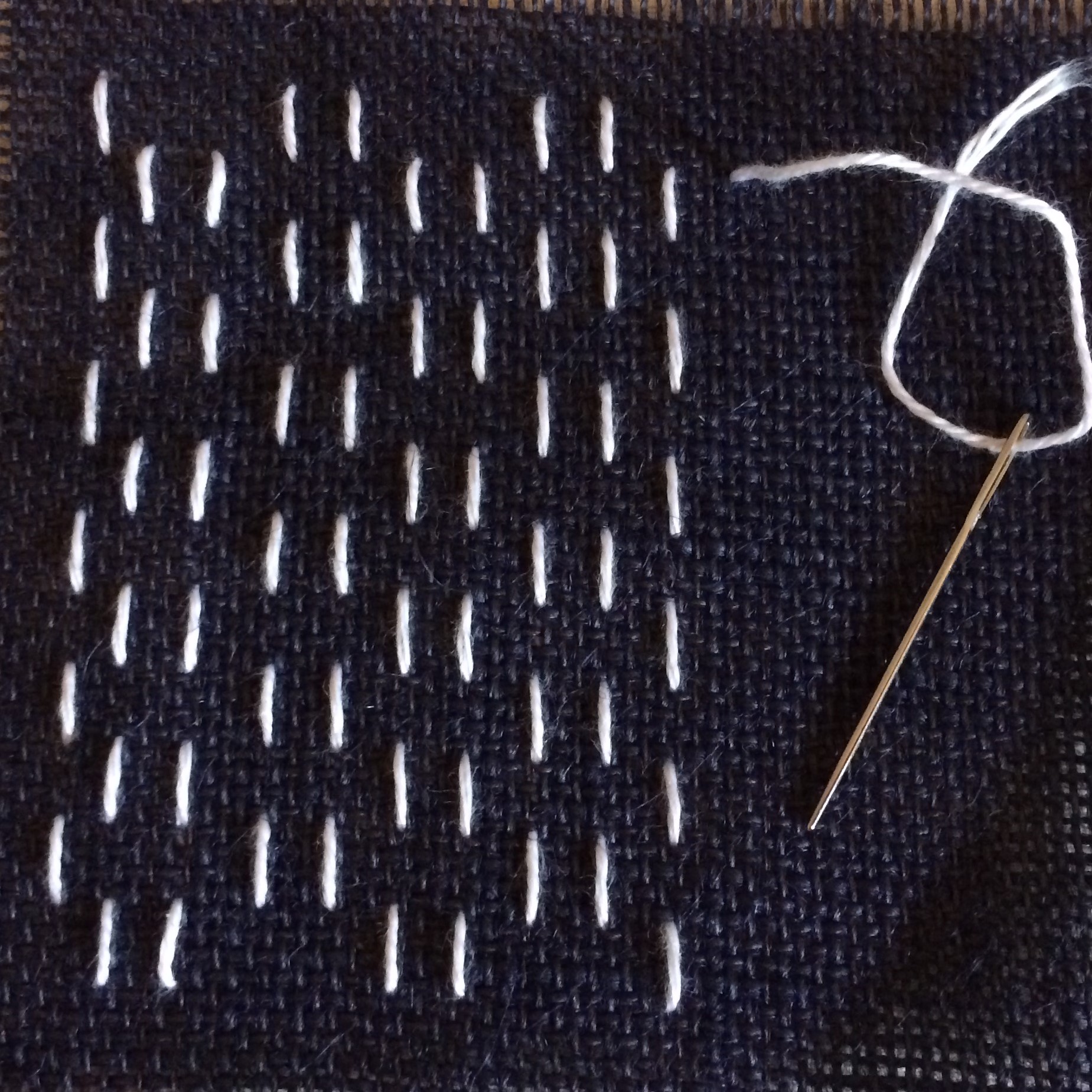}}}
\caption{\label{running} Starting to decorate a small square coaster with hitomezashi. The first ten vertical lines of running stitch have been completed. }
\end{center}
\end{figure}

Coming back to a single line of running stitch, the stitches alternate on the `front' and the `back' of the cloth (see Figure \ref{running}). There are exactly two possible states for the first stitch in a line: either it is present on the front side of the work or it is not present on the front side. It seems natural to encode these two states by 1 (present) and 0 (not). By the nature of running stitch, in each line of stitches the state of that first stitch then determines all the rest. Counting along the line, the odd-numbered stitches are in this same state, and the even-numbered stitches are in the other (opposite) state.

The instructions for working an entire piece of hitomezashi can thus be completely specified by giving two strings of zeroes and ones, the first string giving the state of the first stitch in the verical lines in order and the second doing the same for the horizontal. For example, the piece in Figure \ref{running}, which contains only vertical lines of stitches, can be specified by 0110011001 reading along the bottom, from the left. 

While this simple encoding is not found in the needlework literature, it is so natural that it seems likely that for centuries stitchers have muttered something similar under their breath to keep track of their work as they have stitched. It also enables us to draw on language and notation from the mathematics literature.

For a regular pattern, the sequence of zeroes and ones specifying the vertical lines of stitches can be generated by repeatedly concatenating $v$, a word from the alphabet $\{0,1\}$. We will read the sequence across the bottom of the design, from the left.  In the same way, the horizontal lines of stitches can be specified by repeatedly concatenating a word $w$ from the same binary alphabet, reading up the left-hand side of the design.  The lengths of the words, $|v|$ and $|w|$, need not be the same.

It is useful to define the following: $\overline{w}$ is the word obtained  by the interchange $0\leftrightarrow 1$ in $w$ and $\widetilde{w}$ is the word obtained by reversing the order of the letters in $w$. Both operations are involutions. If $\widetilde{w}=w$, the word $w$ is a palindrome, and if $\widetilde{w}=\overline{w}$, the word $w$ is an antipalindrome.

For example, if $u=01100$, then $\overline{u}=10011$ and $\widetilde{u}=00110$. This word is neither a palindrome nor an antipalindrome and $|u|=5$.


\subsection{Duality}\label{shift}
When hitomezashi patterns are stitched (though not when drawn) a complementary or dual hitomezashi pattern forms on the reverse (Hayes and Seaton, 2020). The gaps in a line of running stitches on one side correspond to stitches on the other, and vice versa.  This may seem unremarkable, but compare it to counted cross-stitch, where the front side of the work consists of cross stitches, but the reverse consists of small vertical stitches. Sashiko is sometimes described as being like blackwork, but there the use of double running stitch results in the same design on both sides of the work (Holden, 2008).

We have found only oblique references to this property in sashiko books. Marquez (2019) remarks that the patterns formed on the back on a piece of hitomezashi are interesting also, and that it `can be reversible with some careful planning'. The near-identical nature of the reverse of some patterns is also mentioned as a distinctive feature of hitomezashi in \textit{Simply Sashiko} (2019/2020), wherein it is explained how decorative kitchen {hanafukin} (flower-cloths)  having two `nice' sides can be made by working through a double layer of cloth, with turning stitches and thread-joins hidden between the two layers. On the other hand, blogger Sawaraka (2017) is delighted that a different (dual) design forms on the reverse when the final diagonal stitches are added to the {kawari hanaj\={u}ji} (variant/changing flower-cross) hitomezashi design.
\begin{figure}
\begin{center}
{\includegraphics[height=5cm]{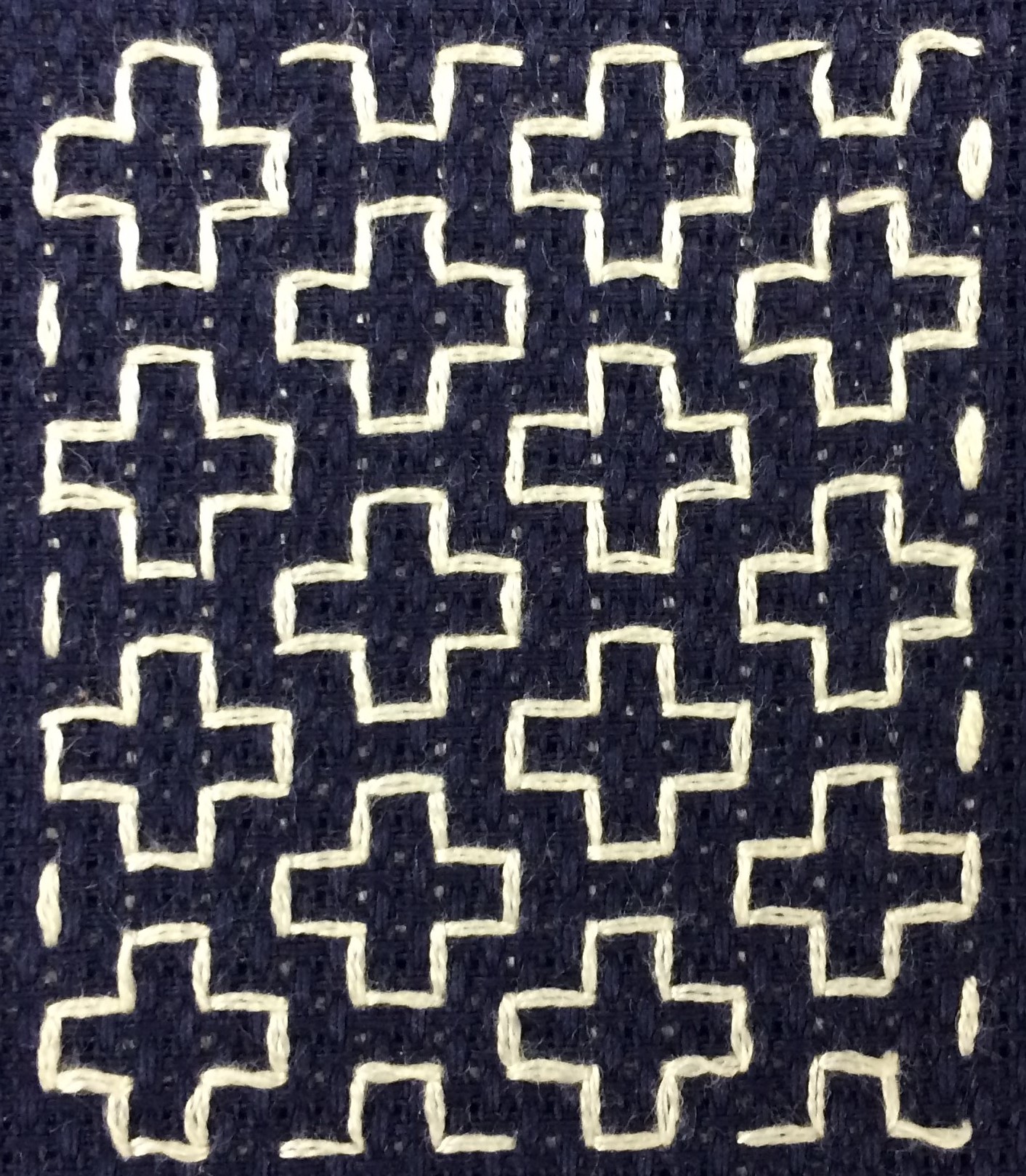}}\quad {\includegraphics[height=5cm]{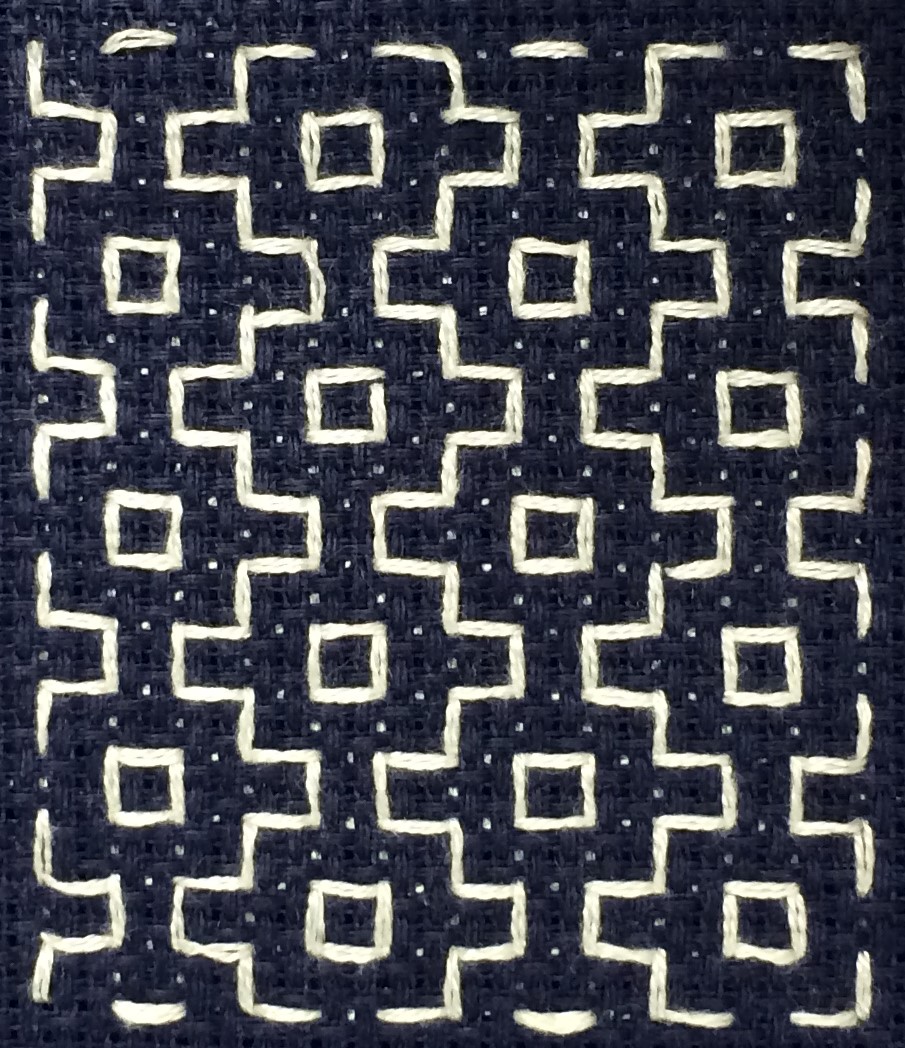}}
\caption{\label{dual} The two sides of a coaster decorated with hitomezashi. On the reverse of a pattern of offset crosses {j\={u}jizashi}, a regular pattern of small squares and stepped lines ({kuchi} and {yamagata}) is formed. This is the dual pattern.}
\end{center}
\end{figure}

 In Figure \ref{dual} the two sides of the same piece of work are shown. When viewing the back of a piece, the way in which it has been turned over can make a difference; this difference will be apparent if there is some asymmetry. The four `gaps' in the centre of each cross in Figure \ref{dual} correspond to the squares on the reverse. To obtain the right-hand image, this piece has been turned over by pivoting on the right or left edge (not by pivoting on the top or bottom edge).

 So that we do not need to consider how a piece of fabric is turned over, we will define the dual pattern to be what would be visible to a viewer, looking over the top of the work into a mirror placed behind the work. Then the stitches behind the left hand side of the work are on the left in the image in the mirror (and so too for right, top and bottom). What is visible is the pattern encoded by $\overline{v}$ and $\overline{w}$.

Hitomezashi designs may be self-dual up to translation. Again, knitters are familiar with this phenomenon. The reverse of a piece of seed stitch is seed stitch, and so too for garter stitch, but the reverse of stocking stitch is purl fabric (i.e., not self-dual). Garter stitch is self-dual up to translation up or down by a row, whereas seed stitch is self-dual up to a shift up or down or left or right by one stitch. When we hereafter describe a pattern as self-dual, we mean `self-dual up to translation'. 

For example, the stitching shown in Figure \ref{running} is self-dual. Recall that it was specified by the word $v=0110011001$. This is a convenient opportunity to introduce the symbol $\epsilon$ for the empty word, with the properties that $\overline{\epsilon}=\epsilon$ and $\widetilde{\epsilon}=\epsilon$, and $|\epsilon|=0$. In Figure \ref{running}, there are no horizontal lines of stitching so that $w=\epsilon$. The dual design is specified by
\[
\overline{v}=\overline{0110011001}=1001100110; \quad \overline{w}=\epsilon
\]
which is a vertical shift by one stitch of the original pattern. This example is simple because $w=\epsilon$. 
 More generally, both the properties of the design’s encoding words under the $0\leftrightarrow 1$ operation and how these interact with each other determine whether a similar or different pattern forms on the reverse.


\subsection{Traditional patterns encoded}\label{trad}
We now give an encoding for some traditional hitomezashi patterns. These encodings are unique (up to cycling and placement of motifs relative to the edge of the worked area). The names of these traditional stitch patterns reflect everyday life and the surroundings in the cold northern Aomori prefecture in the Edo period, the plants, natural phenomena, tools and characters (Mende \& Morishige, 1991).  Many of these are charted in the excellent pattern library of Briscoe (2004).

\begin{figure}
\begin{center}
{\resizebox*{5cm}{!}{\includegraphics{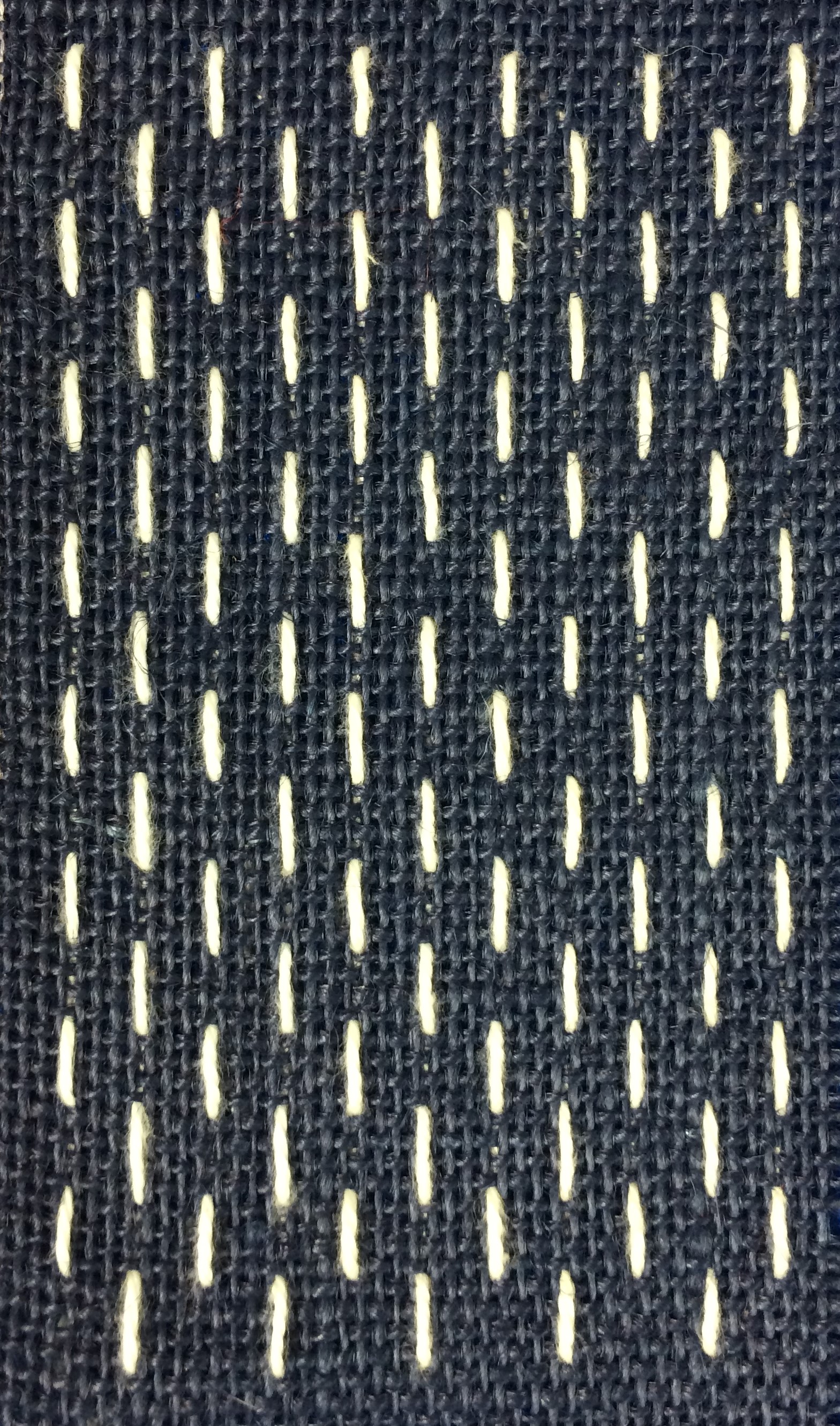}}}
\caption{\label{tate} Vertical lines of stitching obtained by stitching the pattern specified by the word $v=01$ as many times as desired, and to the height desired. This pattern is called {tategushi}. }
\end{center}
\end{figure}

\paragraph{tategushi and yokogushi} These stitch patterns consist of offset vertical (respectively horizontal) lines of stitches, the meaning of their names. {Tategushi} is shown in Figure \ref{tate}. Only one stitch visits each vertex of the grid, and no regions are outlined. For tategushi, $v=01$; $w=\epsilon$. For yokogushi,  $v=\epsilon$; $w=01$.
Each of these designs is self-dual.

\begin{figure}

{\resizebox*{5cm}{!}{\includegraphics{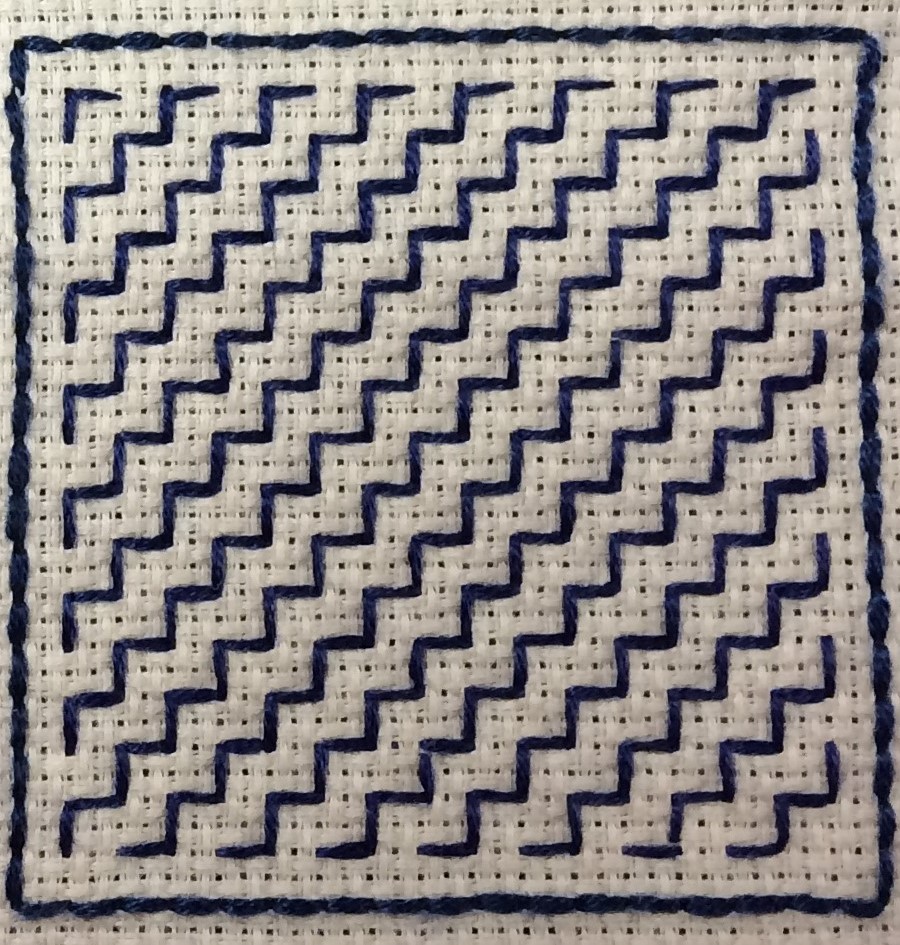}}} \quad \quad \begin{tikzpicture}
  \foreach \x in{0,1,...,4} \foreach \y in {0,1,...,4}{\draw (\x,\y)--({\x+0.5},\y);}

\foreach \x in {0,1,...,4}\foreach \y in {0,1,...,4}{\draw (\x+0.5,{\y+0.5})--({\x+1},{\y+0.5});}
\foreach \y in {0,1,...,4}\foreach \x in {0,1,...,4}{\draw (\x,\y)--(\x, \y+0.5);}
\foreach \y in {0,1,...,4} \foreach \x in {0,1,...,4}{\draw ({\x+0.5},\y+0.5)--(\x+0.5, \y+1);}

\end{tikzpicture}

\hspace{2cm} (a) \hspace{5cm} (b)
\begin{center}
\caption{\label{dan} (a) Linked steps (dan tsunagi) that rise as we move to the right. (b) The same stitch design in the other possible orientation, drawn and to a different scale. This orientation is encoded by $v=10$; $w=10$.}
\end{center}
\end{figure}
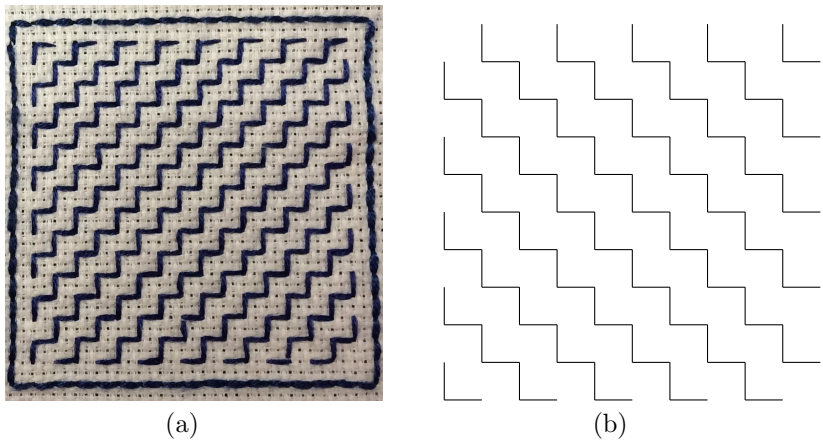

\paragraph{dan tsunagi} The name of this pattern means `linked steps', and it has two possible orientations, as shown in Figures \ref{dan}(a) and \ref{dan}(b).  For steps that rise from southwest to northeast $v=10$; $w=01$.
The dual design is specified by $\overline{v}=01$; $\overline{w}=10$ which is dan tsunagi in the same orientation. That is, the design is self-dual.

\paragraph{kuchizashi} Literally meaning `mouth stitch', the small aligned squares (see the pattern in the top left of Figure \ref{kuchi}) resemble the kanji character mouth (insert unicode character 2F1D here) and also bring to mind the rice and vegetable fields of the T\={o}hoku region. Specified by $v=1$; $w=1$, this design is self-dual.

\begin{figure}
\begin{center}
{{\includegraphics[height=5cm, angle=90]{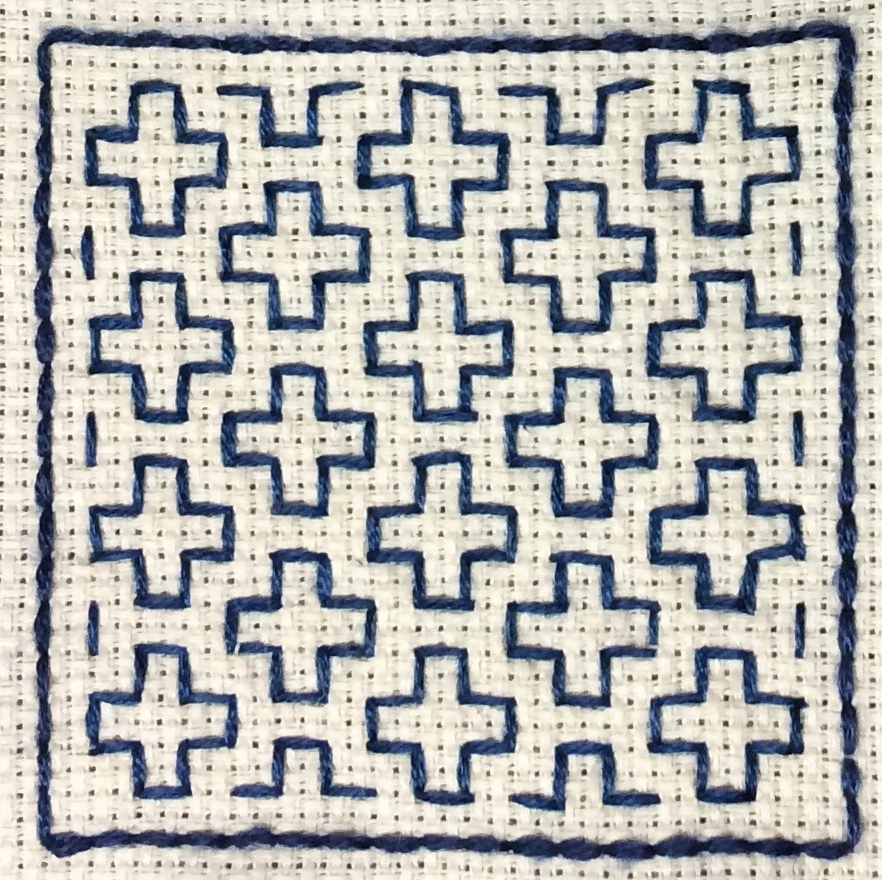}}\qquad {\includegraphics[height=5cm]{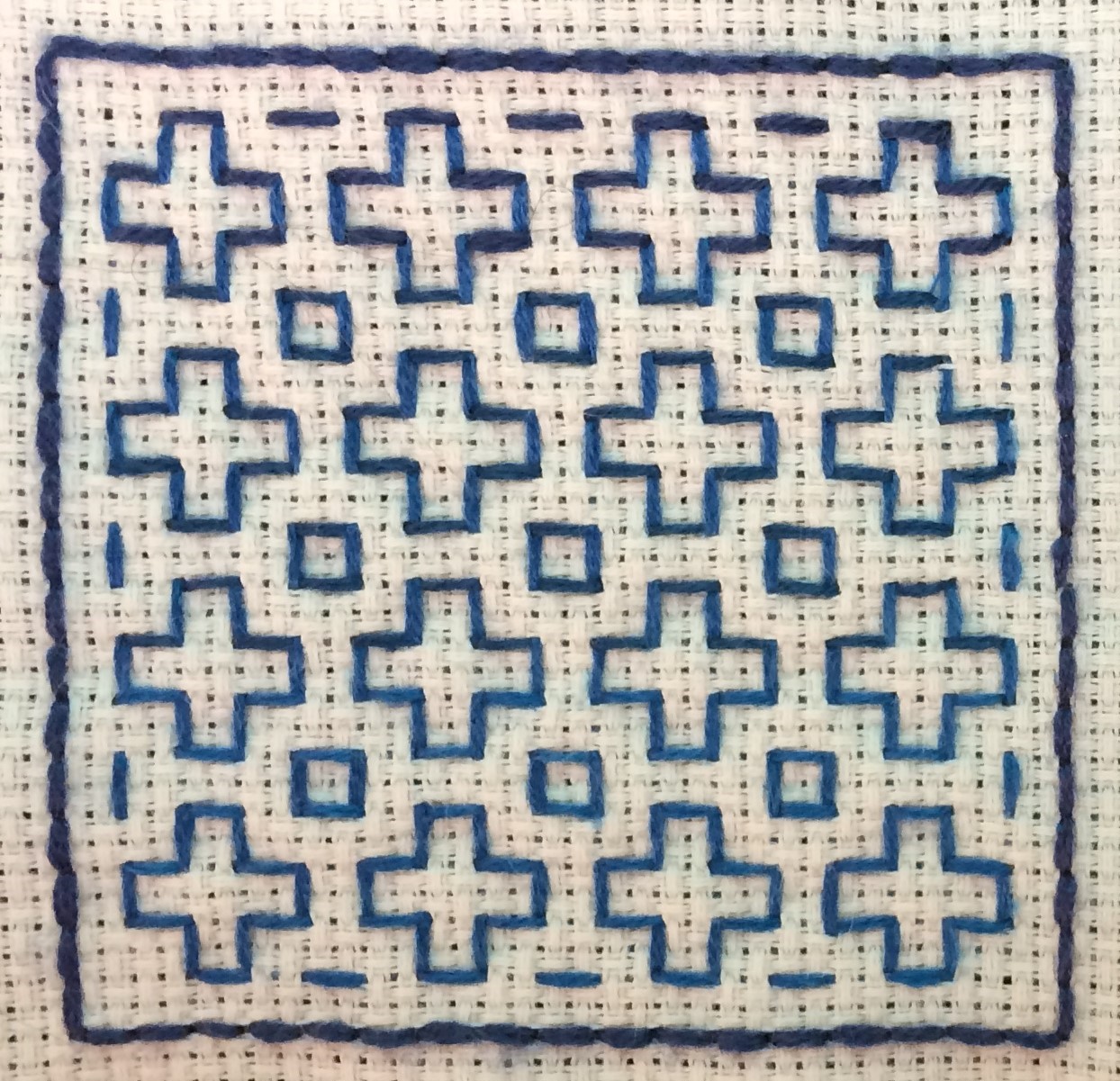}}}

\hspace{1.5cm} (a) \hspace{5cm} (b)
\caption{\label{juji} (a)  Offset ten-crosses ({j\={u}jizashi}). In the orientation shown the encoding is
$v=011$; $w=0110$. (b) Encoded with different words, the crosses are aligned and mouth stitches form between them; this pattern is self-dual. It is discussed in Section 3.}
\end{center}
\end{figure}

\paragraph{j\={u}jizashi}  This pattern (shown in Figure \ref{dual}) gets its name `ten-cross stitch' from the kanji character for the number ten {j\={u}ji}(insert unicode character 2F17 here). We have already observed that it is not self-dual.
The pattern is shown in its other possible orientation (rotated through 90 degrees) in Figure \ref{juji}(a).

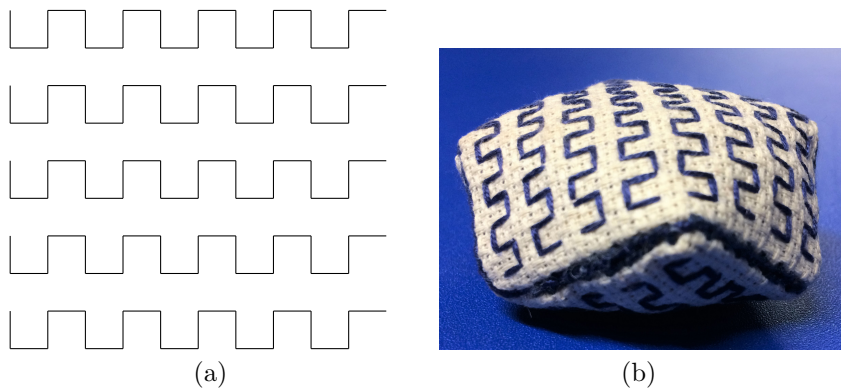
\begin{figure}
\begin{center}
\begin{tikzpicture}
  \foreach \x in{0,1,...,4} \foreach \y in {0,1,...,4}{\draw (\x,\y)--({\x+0.5},\y);}

\foreach \x in {0,1,...,4}\foreach \y in {0,1,...,4}{\draw (\x+0.5,{\y+0.5})--({\x+1},{\y+0.5});}
\foreach \y in {0,1,...,4}\foreach \x in {0,1,...,4}{\draw (\x,\y)--(\x, \y+0.5);}
\foreach \y in {0,1,...,4} \foreach \x in {0,1,...,4}{\draw ({\x+0.5},\y)--(\x+0.5, \y+0.5);}
\end{tikzpicture}\qquad {{\includegraphics[height=4cm]{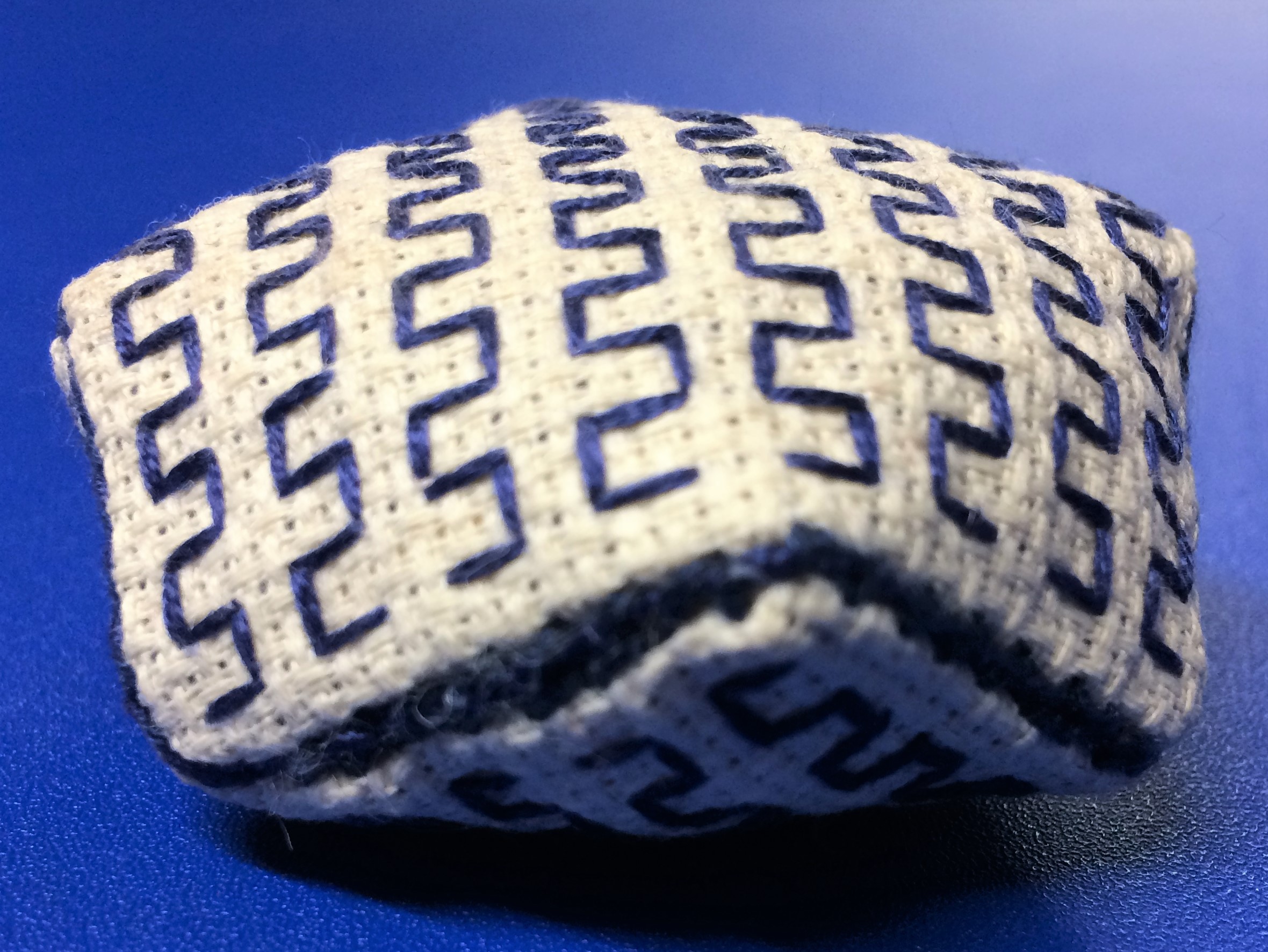}}}
\hspace{1.5cm} (a) \hspace{5cm} (b)
\caption{\label{passes} (a) {hirayama michi} (drawn) with the features aligned, encoded $v=1$; $w=10$. (b) The surface of a pincushion has been stitched with a variant form, wherein the features face away from each other. The encoding is $v=1$; $w=0110$.}
\end{center}
\end{figure}

\paragraph{hirayama michi} The name of this pattern means `passes into the mountain'. This pattern reflects the ups and downs of roads leading from the plains up through the foothills into the high mountains, so the stitcher can choose how many passes to include in their design. When worked as in Figure \ref{passes}(a) the `passes' align and the pattern is self-dual. But a variant of this stitch is to pair two lines in such a way that the `passes' face away from each other as in Figure \ref{passes}(b). It seems opportune to introduce at this point the modifier {kawari} which means `variation'.  The dual pattern in this case is the offset arrangement of squares (mouths) shown in the top right of Figure \ref{kuchi}, {kawari kuchizashi}.

\paragraph{yamagata} Also related to the landscape, this stitch design is called `mountain form'. It resembles the kanji for mountain (insert the unicode character 2F2D).  One way to encode this pattern is as a variant of dan tsunagi, using the word $v=01$ up to the desired position of the mountain, and then using $\widetilde{v}=10$ from that point on, introducing an axis of symmetry. 

\begin{figure}
\begin{center}
{{\includegraphics[height=5cm]{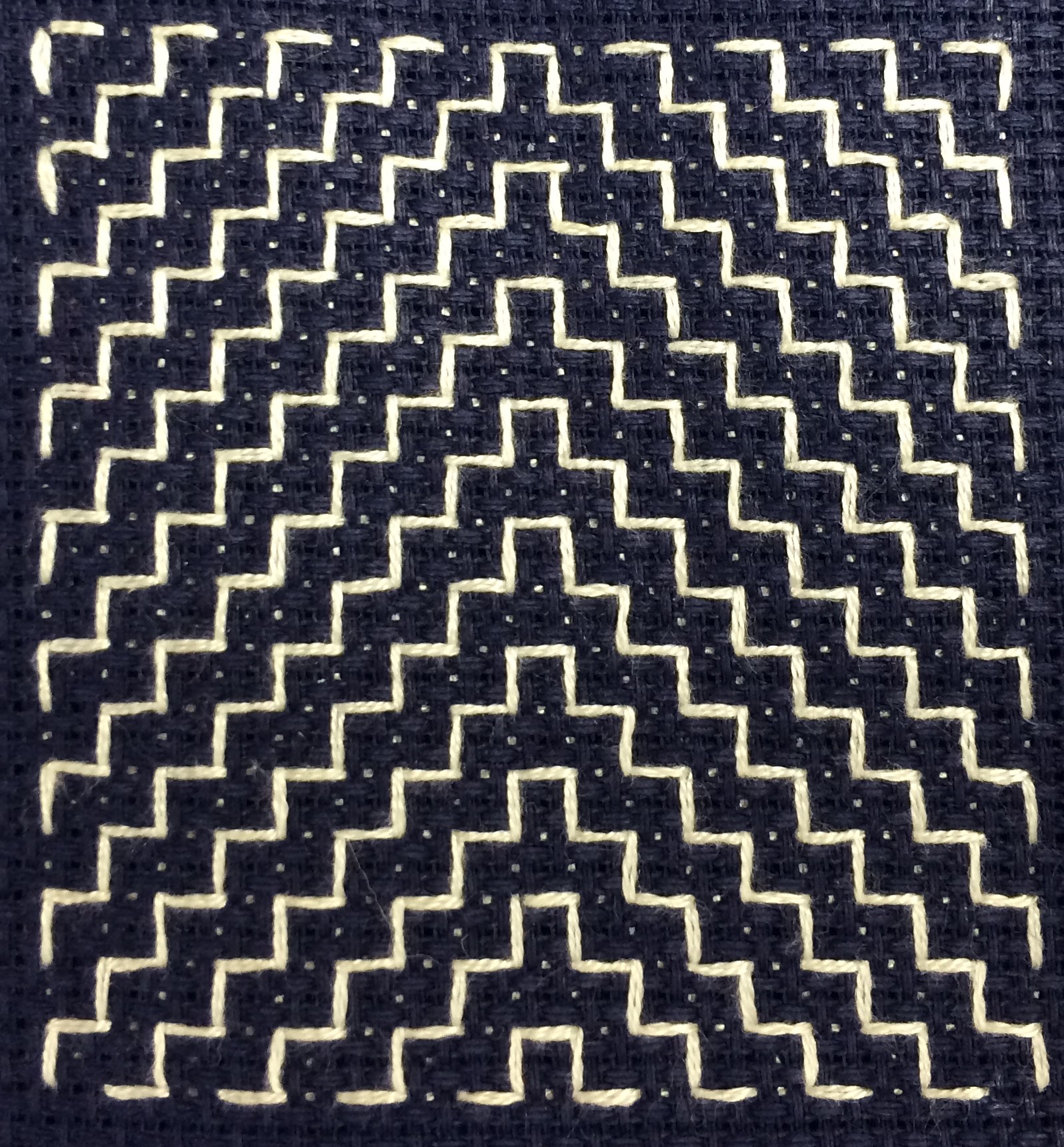}}}\qquad {{\includegraphics[height=5cm]{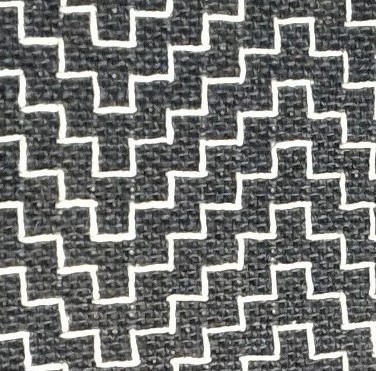}}}

\hspace{1.5cm} (a) \hspace{5cm} (b)
\caption{\label{yama} (a) The basic mountain form {yamagata}, encoded by changing the word from $v=01$ to $\widetilde{v}=10$ at the chosen position of the peak, with $w=01$. (b) Double mountain form {nij\={u} yamagata} which can be worked as a regular repeating design of any height. As shown, the encoding is $v=10101$; $w=10$.  }
\end{center}
\end{figure}

Double mountain form {nij\={u} yamagata} can be of any height; shown in Figure \ref{yama}(b) are mountains of height five stitches (total width ten stitches). Both yamagata and nij\={u} yamagata are self-dual designs. 

\begin{figure}
\begin{center}
{{\includegraphics[height=3.5cm]{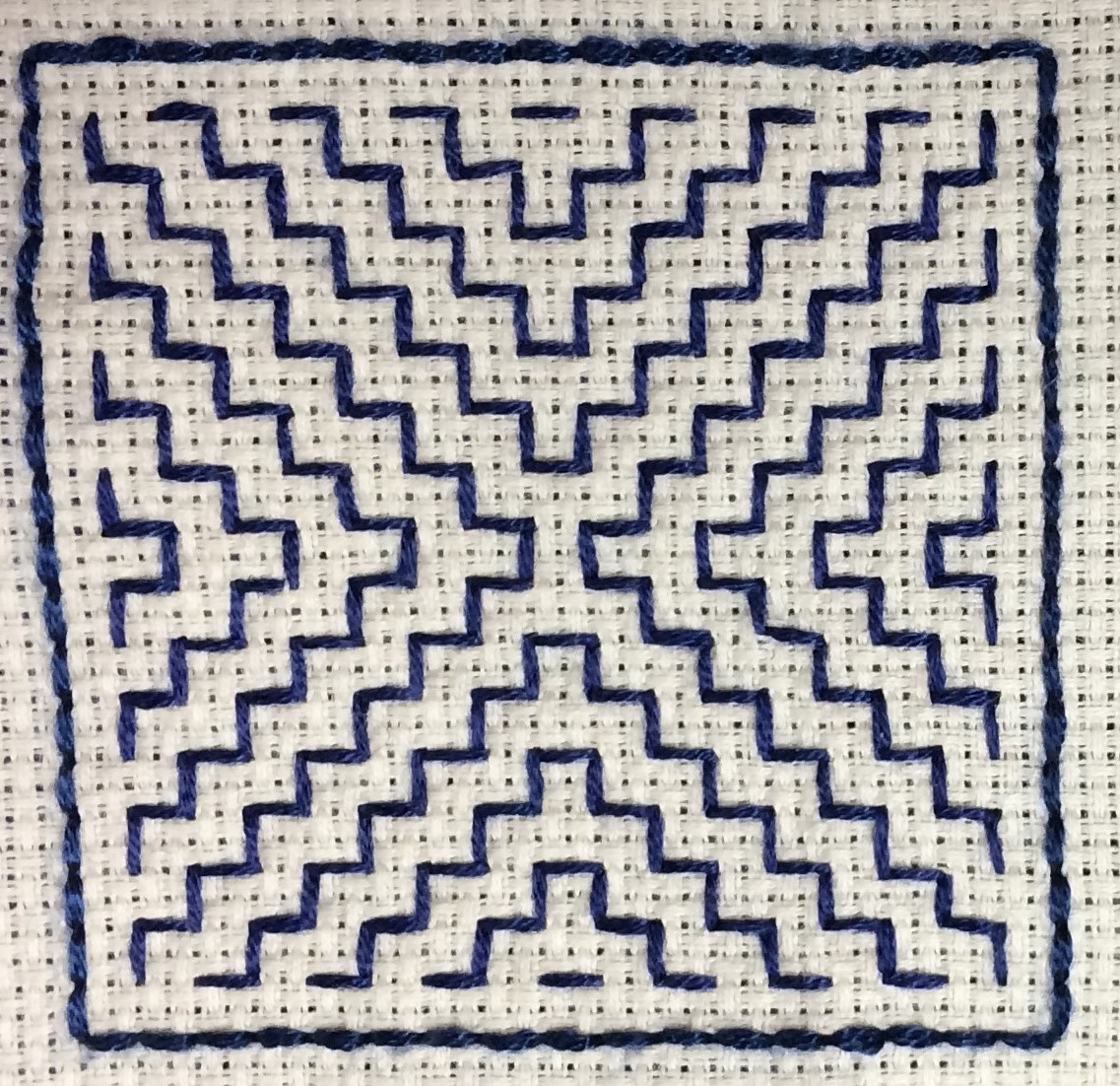}}}\quad {{\includegraphics[height=3.5cm]{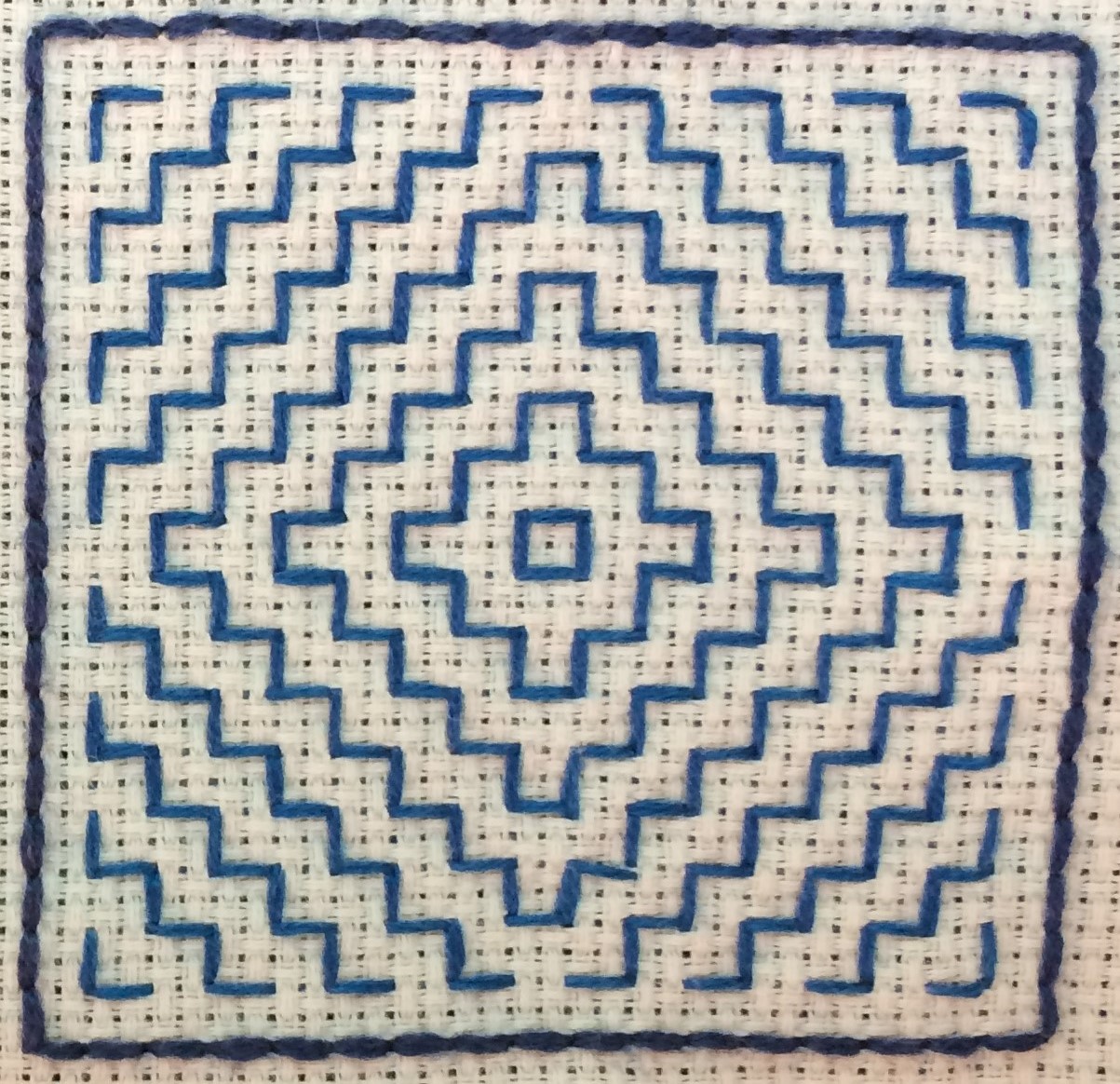}}} \quad {{\includegraphics[height=3.5cm]{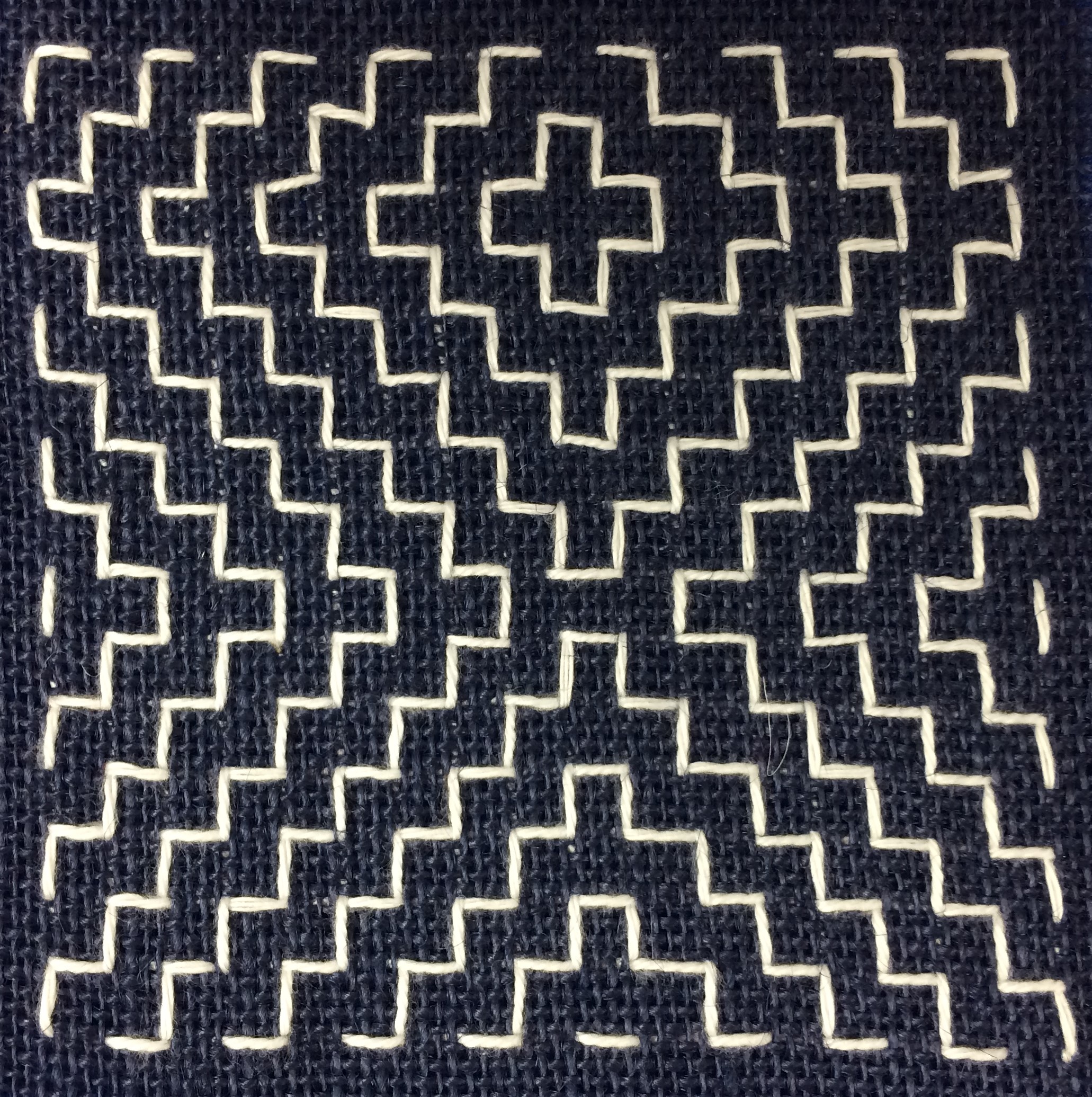}}}

\hspace{1cm} (a) \hspace{3cm} (b) \hspace{3cm}(c)
\caption{\label{yamavar} (a), (b) Further variants of {yamagata}, both extending to the edge of the worked area. (c) The dual of both variants: of (a) if just the lower part is considered, and of (b) if just the upper part is considered.}
\end{center}
\end{figure}

Other variations of {yamagata} are possible, when there is both a horizontal and vertical reflection axis introduced. The reverse of the design in Figure \ref{yamavar}(a) looks almost identical, but careful consideration shows that the dual is a 90 degree rotation of the original (see the lower part of Figure \ref{yamavar}(c)). The dual of the mountain form radiating out from a single square (in Figure \ref{yamavar}(b)) is a mountain form radiating out from a ten-cross (see the upper part of Figure \ref{yamavar}(c)). This second observation is made in \textit{Simply Sashiko} (2019/2020).

\begin{figure}
\begin{center}
{{\includegraphics[height=4.5cm]{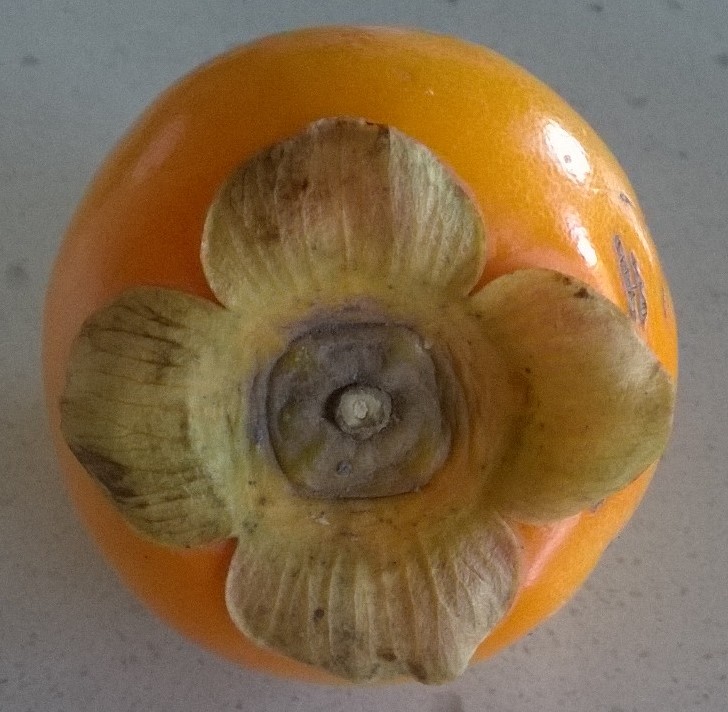}}}\quad {{\includegraphics[height=4.5cm]{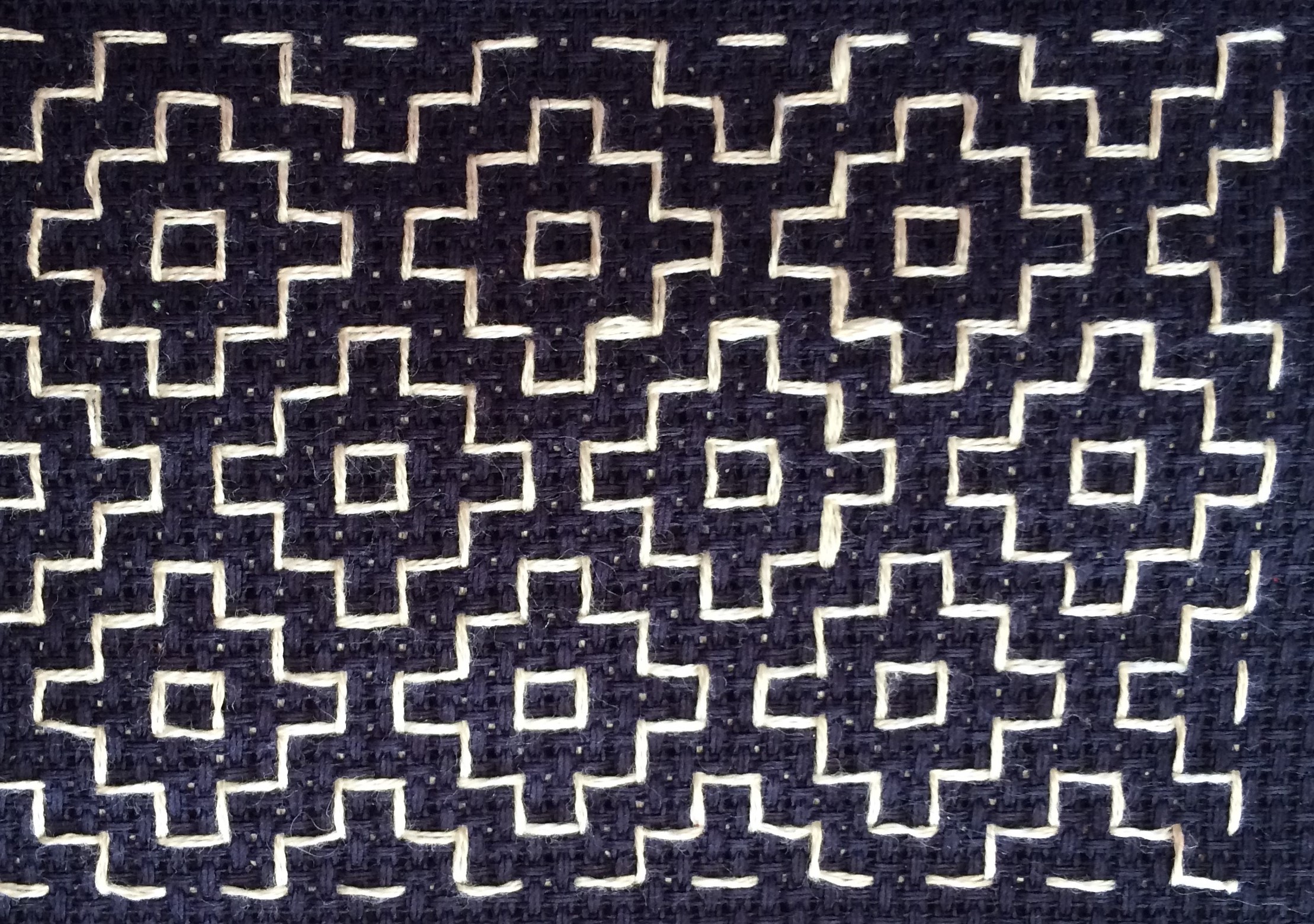}}}
\hspace{1.5cm} (a) \hspace{5cm} (b)

\caption{\label{fruit} (a) A persimmon ({kaki}) fruit with sepals, part of the flower ({hana}) remaining when the fruit forms,  and (b) the stitch named {kakinohanazashi}. The pattern is generated by the words $v=010$ and $w=u\widetilde{u}$ where $u=1010$. }
\end{center}
\end{figure}

\begin{figure}
\begin{center}
{{\includegraphics[height=3cm]{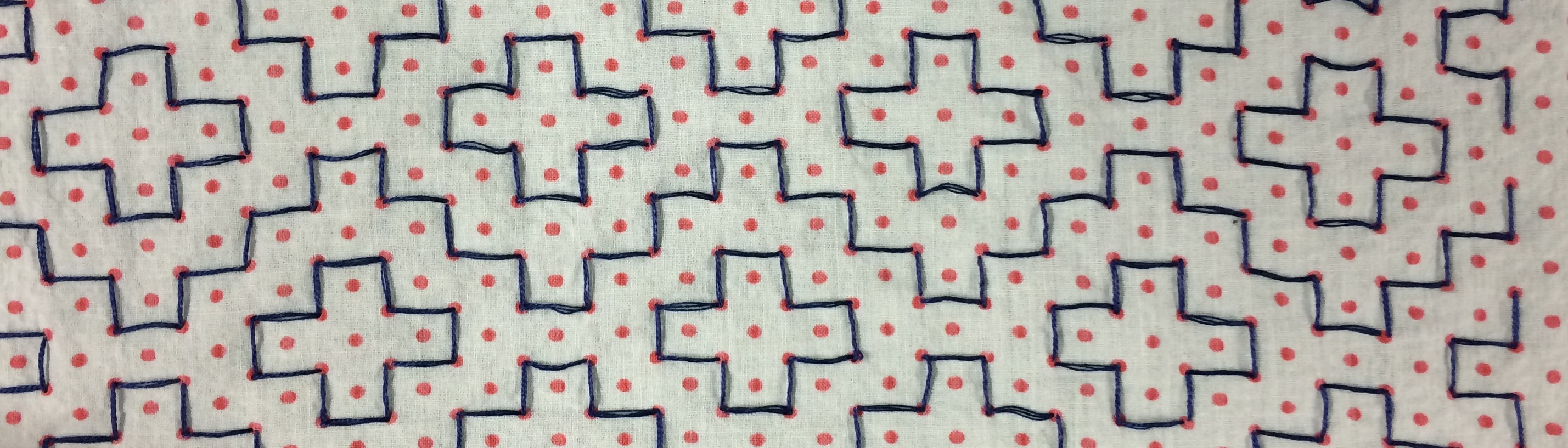}}}
\caption{\label{dualpers} The dual pattern to {kakinohanazashi}. This piece also demonstrates another way to make  hitomezashi stitches of the same length, using a piece of polka-dot or checked fabric to provide a grid.  }
\end{center}
\end{figure}

\paragraph{kakinohanazashi} The oriental persimmon ({kaki}) has an urn-shaped flower ({hana}) with four lobes and four sepals. The sepals remain attached to the attractive golden fruit, which ripens in autumn. Haiku poets sing the praises of this fruit, delicious when eaten at the point of syrupy softness. A repeated design of persimmon flower stitch is shown alongside the fruit in Figure \ref{fruit}; the naming inspiration is obvious. The persimmons are offset, and the pattern is not self-dual. The dual pattern consists of ten-crosses sitting between single lines of nij\={u} yamagata, shown in Figure \ref{dualpers}. 

\begin{figure}
\begin{center}
{{\includegraphics[height=4cm]{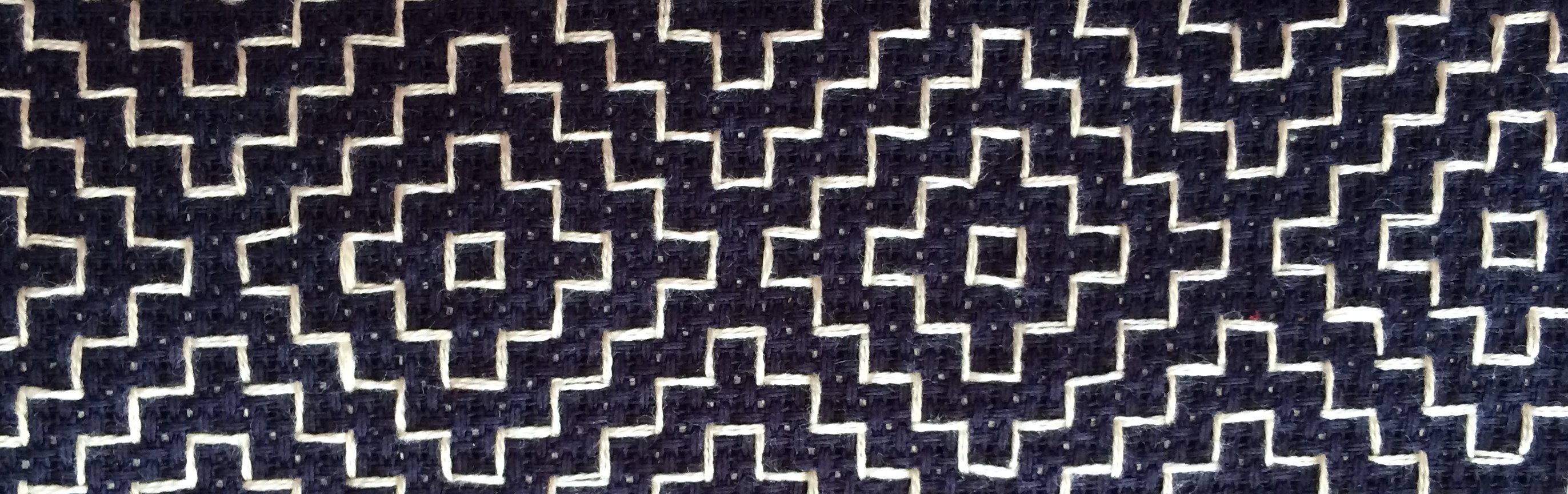}}}
\caption{\label{sanju}  The three closed loops of stitching forming each flower give this stitch pattern the name triple persimmon flower stitch ({sanj\={u} kakinohanazashi}). The encoding is $v=01010$; $w=u\widetilde{u}$ where $u=101010$. }
\end{center}
\end{figure}

\paragraph{sanj\={u} kakinohanazashi} Triple persimmon flower stitch is shown in Figure \ref{sanju}. It is encoded with two more letters in each word, compared to the basic persimmon flower. The dual fabric features a cross at the centre of each motif.

\begin{figure}
\begin{center}
{{\includegraphics[height=3cm]{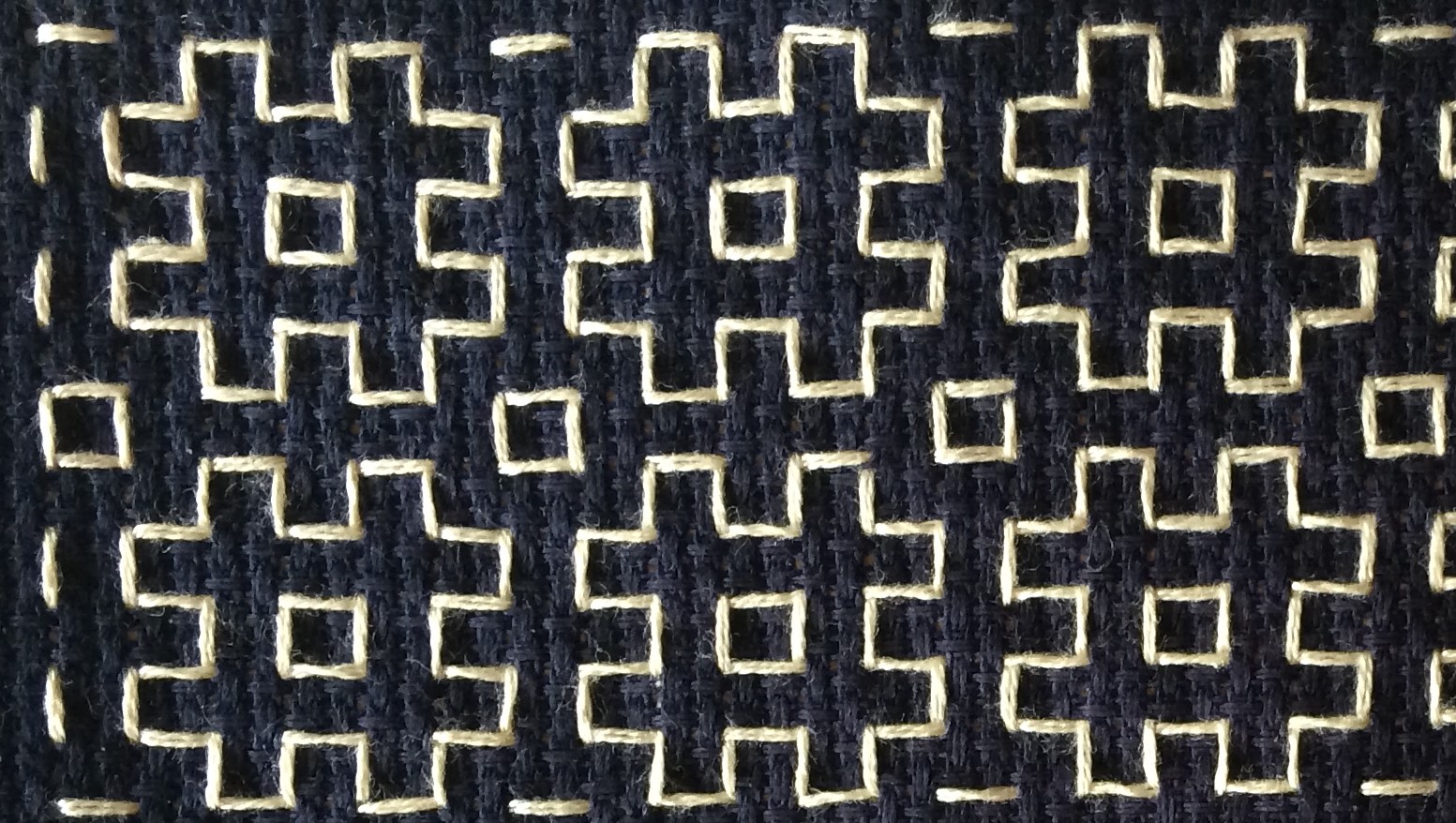}} \quad {\includegraphics[height=3cm]{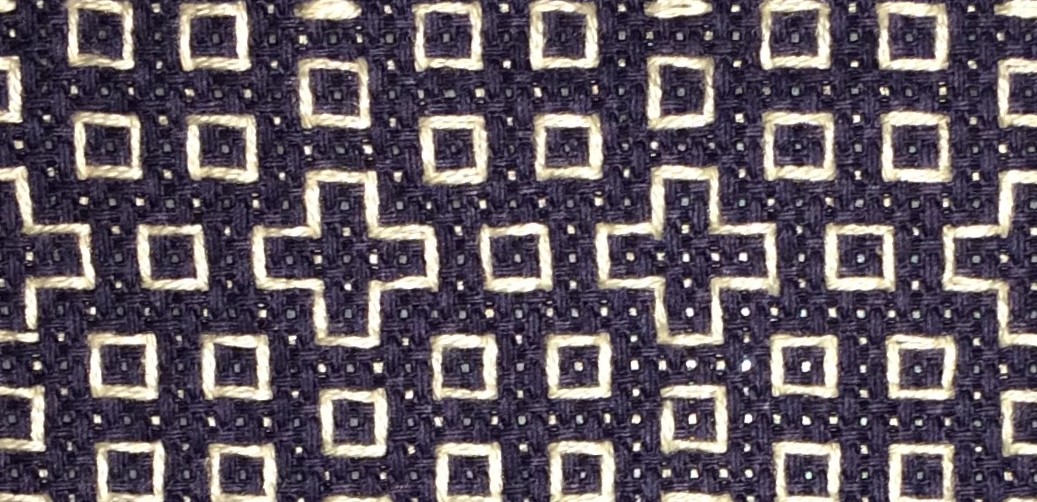}}}
\hspace{1.5cm} (a) \hspace{5cm} (b)
\caption{\label{well} (a) {Igetazashi} or well-kerb stitch.  The encoding is $v=w=u\widetilde{u}$ where $u=100$. (b) The dual fabric. }
\end{center}
\end{figure}

\paragraph{igetazashi} While to modern eyes, this stitch pattern looks like a hashtag, it takes its name from a much older technology. `{i}' means water-well which has kanji (insert unicode character 4E95 here), and {`geta'} refers to the stone-work edging (the kerb/curb). Unlike persimmon flowers, these motifs do not nestle nicely with each other, and  small mouths sit between them. This pattern and its dual are shown in Figure \ref{well}.

\subsection{Loops, lengths and areas} \label{loops}

In some patterns among the selection of traditional designs presented in Section \ref{trad}, closed loops appear, marking out polyominoes (planar figures made by joining identically-sized squares along their edges). In Table 1 their perimeter length, area, width and height are given. (Note, this refers to the size of the largest loop in the stitch pattern, and any enclosed loops are not considered.) The mouth, cross, persimmon and triple persimmon belong to a sequence of polyominoes with area given by the \textit{centred square numbers} 
\[n^2+(n+1)^2=2n(n+1)+1,\ n=0,1,\ldots\]

\begin{table}
    \centering
    \begin{tabular}{c|c|c|c|c}
   Pattern name& Perimeter&Area&Height&Width\\
    \hline
 kuchizashi      & 4 &1&1&1\\
  j\={u}jizashi      &12 &5&3&3\\
kakinohanazashi   &20&13&5&5  \\
dual sanj\={u} kakinohanazashi&28&25&7&7\\
sanj\={u} kakinohanazashi &36&41&9&9\\\hline
igetazashi&28&17&5&5
\end{tabular}
    \caption{Features of the largest loops appearing in some traditional hitomezashi patterns}
    \label{features}
\end{table}
     
Pete (2008) showed (in the context of corner percolation) that the width and height of any hitomezashi loop will be odd. Since there must be one more boundary stitch than there are squares, both of the subwords encoding the part of the pattern containing the loop are of even length.
Defant and Kravitz (2022) have further shown that for an arbitrary hitomezashi loop:
 \begin{itemize}
     \item The enclosed area is congruent to 1 modulo 4;
     \item The length of the perimeter is congruent to 4 modulo 8.
   \end{itemize}
These features are, of course, displayed by the patterns listed in Table \ref{features}.


%

\section{Fibonacci snowflakes and Pell persimmons}
\subsection{Fibonacci snowflakes are hitomezashi loops}
\begin{figure}
\begin{center}
{{\includegraphics[height=8cm]{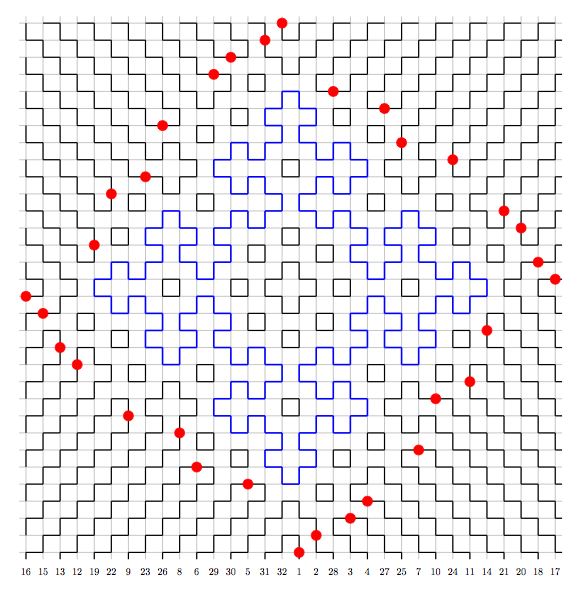}}}
\caption{\label{labbe} Within this fully packed loop diagram for the permutation  
indicated by the numbers along the bottom, 
Fibonacci snowflakes of orders one, two and three occur. Reproduced from Labb\'{e} (2010) with permission.}
\end{center}
\end{figure}

Investigation of a possible connection between hitomezashi loops and the fully packed loop models of statistical mechanics led us to the remarkable image found in the blog of Labb\'{e} (2010) reproduced as Figure \ref{labbe}. Mouths and crosses sit inside a larger polyomino, and in such a way that the central part of the design can be stitched using hitomezashi running stitches (Hayes \& Seaton, 2020). This larger polyomino is one of the Fibonacci snowflakes. 

\begin{figure}
\begin{center}
\begin{tikzpicture}
\draw(0,0)--(0,0.5);
\draw(0,0)--(0.5,0);
\draw(0.5,0)--(0.5,0.5);
\draw(0,0.5)--(0.5,0.5);
\draw(2,0)--(2,0.5);
\draw(2,0)--(2.5,0);
\draw(2.5,0)--(2.5,0.5);
\draw(1.5,0.5)--(2,0.5);
\draw(1.5,0.5)--(1.5,1);
\draw(2,1)--(2,1.5);
\draw(1.5,1)--(2,1);
\draw(2,1.5)--(2.5,1.5);
\draw(2.5,0.5)--(3,0.5);
\draw(3,0.5)--(3,1);
\draw(2.5,1)--(3,1);
\draw(2.5,1)--(2.5,1.5);
\draw(5,0)--(5,0.5);
\draw(5,0)--(5.5,0);
\draw(5.5,0)--(5.5,0.5);
\draw(4.5,0.5)--(4.5,1);
\draw(6,0.5)--(6,1);
\draw(4.5,1.5)--(4.5,2);
\draw(5,1)--(5,1.5);
\draw(5.5,1)--(5.5,1.5);
\draw(6,1.5)--(6,2);
\draw(4,2)--(4.5,2);
 \foreach \x in{4,5} \foreach \y in {0.5,1, 1.5,3,3.5,4}{\draw (\x+0.5,\y)--({\x+1},\y);}
\foreach \x in {3,7.5}\foreach \y in {2}{\draw (\x,{\y})--({\x},{\y+0.5});}
\foreach \x in {3.5,6.5}\foreach \y in {1.5,3}{\draw (\x,{\y})--({\x+0.5},{\y});}
\foreach \x in {3.5,4,6.5,7}{\draw (\x,1.5)--({\x},2);} 
\foreach \x in {3.5,4,6.5,7}{\draw (\x,2.5)--({\x},3);}
\foreach \x in {4.5,6}{\draw (\x,2.5)--({\x},3);}
\foreach \x in {5,5.5}{\draw (\x,3)--({\x},3.5);} 
\foreach \x in {5,5.5}{\draw (\x,4)--({\x},4.5);}
\foreach \y in {2,2.5}\foreach \x in {3,4,6,7}{\draw (\x,\y)--(\x+0.5, \y);}
\foreach \y in {3.5}\foreach \x in {4.5,6}{\draw (\x,\y)--(\x, \y+0.5);}
\draw(5,4.5)--(5.5,4.5);
\draw[dotted](5,2)--(5.5,2)--(5.5,2.5)--(5,2.5)--(5,2);

\end{tikzpicture}

\caption{\label{snow} The first three Fibonacci snowflakes are shown. Their outline is defined using words generated from equation (\ref{fibword}). The dashed line indicates how the order one object forms inside the order three object if it is, rather, stitched using vertical and horizontal lines of running stitch.}
\end{center}
\end{figure}
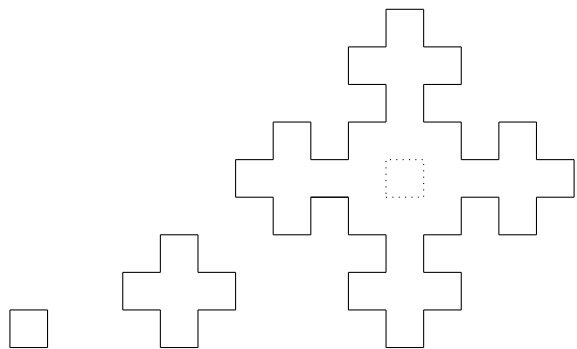

  \textit{Fibonacci snowflakes}, also sometimes called \textit{Fibonacci tiles} or \textit{Fibonacci polyominoes}, were introduced by Blondin-Mass\'{e} et al. (2009). These are polyominoes defined using a subset of a family of words on two letters concatenated in such a way that the number of letters in the word $q_n$ is the $n$-th \textit{Fibonacci number}, $F_n$. (These words are related to, but not identical to, the sequence of finite words leading in the limit to the Fibonacci word.) The alphabet used is $\{L, R\}$ indicating unit length steps to the left or right on a square grid.  In particular
 \begin{equation}
 q_n=\begin{cases}
  q_{n-1}q_{n-2}& n \equiv 2 \mod{3}\\
 q_{n-1}\overline{q_{n-2}}& n \equiv 0,1 \mod{3} \label{fibword}
\end{cases}
 \end{equation}
where $\overline{q}$ denotes $L\leftrightarrow R$ in $q$, and $q_0=\epsilon$ and $q_1=R$. The $n$-th Fibonacci snowflake has as its boundary curve the four-fold concatenation $(q_{3n+1})^4$ (Blondin-Mass\'{e}, Brlek, Garon, et al., 2011).  Note that these steps are drawn without lifting one's pen from the page; they are not created in two stages as hitomezashi loops are. The first few are shown in Figure \ref{snow}. 

These objects were identified in a systematic search for polyominoes that tile the plane (Blondin-Mass\'{e} et al., 2009).  They have a number of interesting properties. The Fibonacci snowflakes are pseudo-squares, tiling the plane in such a way that each polyonimo is bordered by four copies of itself, as indicated in Figure \ref{tess}. The perimeter length of the $n$-th snowflake is $4F_{3n-2}$. 

From the recursion relation
\[
F_n=F_{n-1}+F_{n-2}; \quad F_0=0, \ F_1=1
\] it is established that
\[
F_{3k-2}=3F_{3k-5}+2F_{3k-6}
\]
so that (by induction), $F_{3n-2}$ is odd.

The area of each Fibonacci snowflake is given by an odd-index Pell number (Blondin-Mass\'{e}, Brlek, Labb\'{e}, et al., 2011).  The \textit{Pell numbers} are defined by
\begin{equation}
 P_n=2P_{n-1}+P_{n-2}; \quad P_0=0, P_1=1.   \label{pell}
\end{equation}
The first few are:
\[
0,\ 1,\ 2,\ 5,\ 12,\ 29,\ 70, \ldots
\]
By induction, since
\[
P_{2k+1}=4(P_{2k-1}+P_{2k-2})+P_{2k-3},
\]
the odd-index Pell numbers are congruent to 1 modulo 4.

Thus the length and area of the Fibonacci snowflakes satisfy the conditions shown by Defant and Kravitz (2022) for hitomezashi loops (see Section \ref{loops}).

The argument leading to the area formula presented by Blondin-Mass\'{e}, Brlek, Labb\'{e} et al. (2011) involves all Pell numbers, not only those with odd index. We make the related observation that the width (in boundary stitches) of the $n$-th Fibonacci snowflake (considered as a hitomezashi loop) is twice the Pell number $P_n$.

The naming as `snowflake' appears inspired by similarity to the Koch snowflake; the four-fold symmetry is not that of a natural snowflake (which is six-fold). Blondin-Mass\'{e} et al. (2012)  showed that the fractal obtained from the series of Fibonacci snowflakes as  $n$ gets large has fractal dimension 
\[\frac{\log(2+\sqrt{5})}{\log(1+\sqrt{2})}.\]

 Among the psuedo-square polyominoes of the paper by Blondin-Mass\'{e} et al. (2010) and also those in the work of Ramirez et al. (2014), we have identified others which can be stitched in hitomezashi, but we focus here only on the Fibonacci snowflakes.

\subsection{From loops to overall designs}
\begin{figure}
\begin{center}
\begin{tikzpicture}

\draw[gray](6,2.5)--(6,3);
\draw[gray](6,2.5)--(6.5,2.5);
\draw[gray](6.5,2.5)--(6.5,3);
\draw[gray](5.5,3)--(5.5,3.5);
\draw[gray](7,3)--(7,3.5);
\draw[gray](5.5,4)--(5.5,4.5);
\draw[gray](6,3.5)--(6,4);
\draw[gray](6.5,3.5)--(6.5,4);
\draw[gray](7,4)--(7,4.5);
\draw[gray](5,4.5)--(5.5,4.5);
 \foreach \x in{5,6} \foreach \y in {3,3.5, 4,5.5,6,6.5}{\draw[gray] (\x+0.5,\y)--({\x+1},\y);}
\foreach \x in {4,8.5}\foreach \y in {4.5}{\draw[gray] (\x,{\y})--({\x},{\y+0.5});}
\foreach \x in {4.5,7.5}\foreach \y in {4,5.5}{\draw[gray] (\x,{\y})--({\x+0.5},{\y});}
\foreach \x in {4.5,5,7.5,8}{\draw[gray] (\x,4)--({\x},4.5);} 
\foreach \x in {4.5,5,7.5,8}{\draw[gray] (\x,5)--({\x},5.5);}
\foreach \x in {5.5,7}{\draw[gray] (\x,5)--({\x},5.5);}
\foreach \x in {6,6.5}{\draw[gray] (\x,5.5)--({\x},6);} 
\foreach \x in {6,6.5}{\draw (\x,6.5)--({\x},7);}
\foreach \y in {4.5,5}\foreach \x in {4,5,7,8}{\draw[gray] (\x,\y)--(\x+0.5, \y);}
\foreach \y in {6}\foreach \x in {5.5,7}{\draw[gray] (\x,\y)--(\x, \y+0.5);}
\draw[gray](6,7)--(6.5,7);
%
\draw[gray](2.5,1)--(2.5,1.5);
\draw[gray](2.5,1)--(3,1);
\draw[gray](3,1)--(3,1.5);
\draw[gray](2,1.5)--(2,2);
\draw[gray](3.5,1.5)--(3.5,2);
\draw[gray](2,2.5)--(2,3);
\draw[gray](2.5,2)--(2.5,2.5);
\draw[gray](3,2)--(3,2.5);
\draw[gray](3.5,2.5)--(3.5,3);
\draw[gray](1.5,3)--(2,3);
 \foreach \x in{1.5,2.5} \foreach \y in {1.5,2, 2.5,4,4.5,5}{\draw[gray] (\x+0.5,\y)--({\x+1},\y);}
\foreach \x in {0.5,5}\foreach \y in {3}{\draw[gray] (\x,{\y})--({\x},{\y+0.5});}
\foreach \x in {1,4}\foreach \y in {2.5,4}{\draw[gray] (\x,{\y})--({\x+0.5},{\y});}
\foreach \x in {1,1.5,4,4.5}{\draw[gray] (\x,2.5)--({\x},3);} 
\foreach \x in {1,1.5,4,4.5}{\draw[gray] (\x,3.5)--({\x},4);}
\foreach \x in {2,3.5}{\draw[gray] (\x,3.5)--({\x},4);}
\foreach \x in {2.5,3}{\draw[gray] (\x,4)--({\x},4.5);} 
\foreach \x in {2.5,3}{\draw[gray] (\x,5)--({\x},5.5);}
\foreach \y in {3,3.5}\foreach \x in {0.5,1.5,3.5,4.5}{\draw[gray] (\x,\y)--(\x+0.5, \y);}
\foreach \y in {4.5}\foreach \x in {2,3.5}{\draw[gray] (\x,\y)--(\x, \y+0.5);}
\draw[gray](2.5,5.5)--(3,5.5);
\draw(5,0)--(5,0.5);
\draw(5,0)--(5.5,0);
\draw(5.5,0)--(5.5,0.5);
\draw(4.5,0.5)--(4.5,1);
\draw(6,0.5)--(6,1);
\draw(4.5,1.5)--(4.5,2);
\draw(5,1)--(5,1.5);
\draw(5.5,1)--(5.5,1.5);
\draw(6,1.5)--(6,2);
\draw(4,2)--(4.5,2);
 \foreach \x in{4,5} \foreach \y in {0.5,1, 1.5,3,3.5,4}{\draw (\x+0.5,\y)--({\x+1},\y);}
\foreach \x in {3,7.5}\foreach \y in {2}{\draw (\x,{\y})--({\x},{\y+0.5});}
\foreach \x in {3.5,6.5}\foreach \y in {1.5,3}{\draw (\x,{\y})--({\x+0.5},{\y});}
\foreach \x in {3.5,4,6.5,7}{\draw (\x,1.5)--({\x},2);} 
\foreach \x in {3.5,4,6.5,7}{\draw (\x,2.5)--({\x},3);}
\foreach \x in {4.5,6}{\draw (\x,2.5)--({\x},3);}
\foreach \x in {5,5.5}{\draw (\x,3)--({\x},3.5);} 
\foreach \x in {5,5.5}{\draw (\x,4)--({\x},4.5);}
\foreach \y in {2,2.5}\foreach \x in {3,4,6,7}{\draw (\x,\y)--(\x+0.5, \y);}
\foreach \y in {3.5}\foreach \x in {4.5,6}{\draw (\x,\y)--(\x, \y+0.5);}
\draw(5,4.5)--(5.5,4.5);
\draw[gray](4,-2.5)--(4,-2);
\draw[gray](4,-2.5)--(4.5,-2.5);
\draw[gray](4.5,-2.5)--(4.5,-2);
\draw[gray](3.5,-2)--(3.5,-1.5);
\draw[gray](5,-2)--(5,-1.5);
\draw[gray](3.5,-1)--(3.5,-0.5);
\draw[gray](4,-1.5)--(4,-1);
\draw[gray](4.5,-1.5)--(4.5,-1);
\draw[gray](5,-1)--(5,-0.5);
\draw[gray](3,-0.5)--(3.5,-0.5);
 \foreach \x in{3,4} \foreach \y in {-2,-1.5, -1,0.5,1,1.5}{\draw[gray] (\x+0.5,\y)--({\x+1},\y);}
\foreach \x in {2,6.5}\foreach \y in {-0.5}{\draw[gray] (\x,{\y})--({\x},{\y+0.5});}
\foreach \x in {2.5,5.5}\foreach \y in {-1,0.5}{\draw[gray] (\x,{\y})--({\x+0.5},{\y});}
\foreach \x in {2.5,3,5.5,6}{\draw[gray] (\x,-1)--({\x},-0.5);} 
\foreach \x in {2.5,3,5.5,6}{\draw[gray] (\x,0)--({\x},0.5);}
\foreach \x in {3.5,5}{\draw[gray] (\x,0)--({\x},0.5);}
\foreach \x in {4,4.5}{\draw[gray] (\x,0.5)--({\x},1);} 
\foreach \x in {4,4.5}{\draw[gray] (\x,1.5)--({\x},2);}
\foreach \y in {-0.5,0}\foreach \x in {2,3,5,6}{\draw[gray] (\x,\y)--(\x+0.5, \y);}
\foreach \y in {1}\foreach \x in {3.5,5}{\draw[gray] (\x,\y)--(\x, \y+0.5);}
\draw[gray](4,2)--(4.5,2);
%
\draw[gray](7.5,-1)--(7.5,-0.5);
\draw[gray](7.5,-1)--(8,-1);
\draw[gray](8,-1)--(8,-0.5);
\draw[gray](7,-0.5)--(7,0);
\draw[gray](8.5,-0.5)--(8.5,0);
\draw[gray](7,0.5)--(7,1);
\draw[gray](7.5,0)--(7.5,0.5);
\draw[gray](8,0)--(8,0.5);
\draw[gray](8.5,0.5)--(8.5,1);
\draw[gray](6.5,1)--(7,1);
 \foreach \x in{6.5,7.5} \foreach \y in {-0.5,0, 0.5,2,2.5,3}{\draw[gray] (\x+0.5,\y)--({\x+1},\y);}
\foreach \x in {5.5,10}\foreach \y in {1}{\draw[gray] (\x,{\y})--({\x},{\y+0.5});}
\foreach \x in {6,9}\foreach \y in {0.5,2}{\draw[gray] (\x,{\y})--({\x+0.5},{\y});}
\foreach \x in {6,6.5,9,9.5}{\draw[gray] (\x,0.5)--({\x},1);} 
\foreach \x in {6,6.5,9,9.5}{\draw[gray] (\x,1.5)--({\x},2);}
\foreach \x in {7,8.5}{\draw[gray] (\x,1.5)--({\x},2);}
\foreach \x in {7.5,8}{\draw[gray] (\x,2)--({\x},2.5);} 
\foreach \x in {7.5,8}{\draw[gray] (\x,3)--({\x},3.5);}
\foreach \y in {1,1.5}\foreach \x in {5.5,6.5,8.5,9.5}{\draw[gray] (\x,\y)--(\x+0.5, \y);}
\foreach \y in {2.5}\foreach \x in {7,8.5}{\draw[gray] (\x,\y)--(\x, \y+0.5);}
\draw[gray](7.5,3.5)--(8,3.5);
\end{tikzpicture}
\caption{\label{tess} Tiling of the plane by order 3 Fibonacci snowflakes. Note however that this design cannot be stitched as hitomezashi, due to the vertices of order four. This figure is inspired by Figure 2 of the paper by Blondin-Mass\'{e}, Brlek, Labb\'{e} et al. (2011).}
\end{center}
\end{figure}
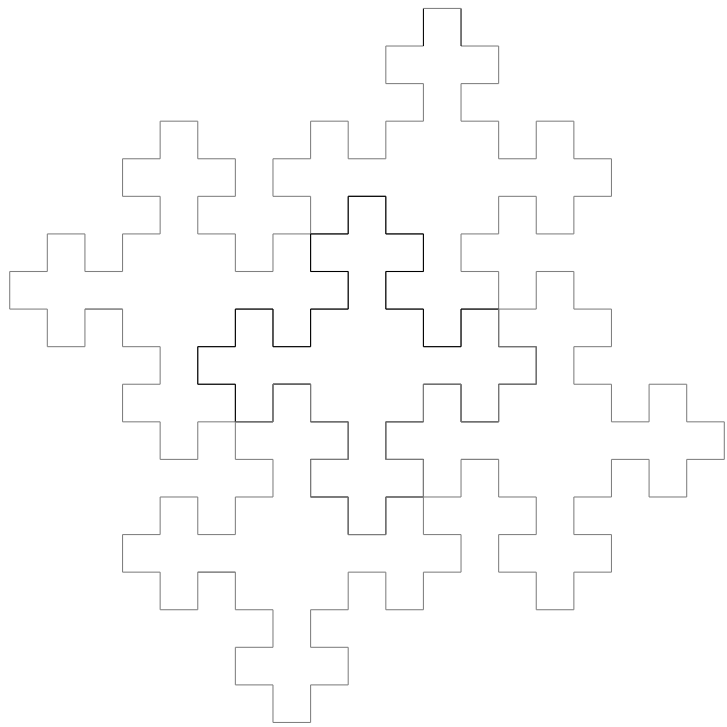

When Fibonacci snowflakes tile the plane as shown in Figure \ref{tess}, their boundaries touch and meet in ways that are not compatible with hitomezashi stitching rules. That is, an overall hitomezashi design featuring Fibonacci snowflakes cannot be based on these tessellations. When one attempts to stitch a single Fibonacci snowflake using lines of running stitch, a flurry of smaller snowflakes appear within and around it, limiting how closely they can be packed together.

\begin{figure}
\begin{center}
{{\includegraphics[height=6cm]{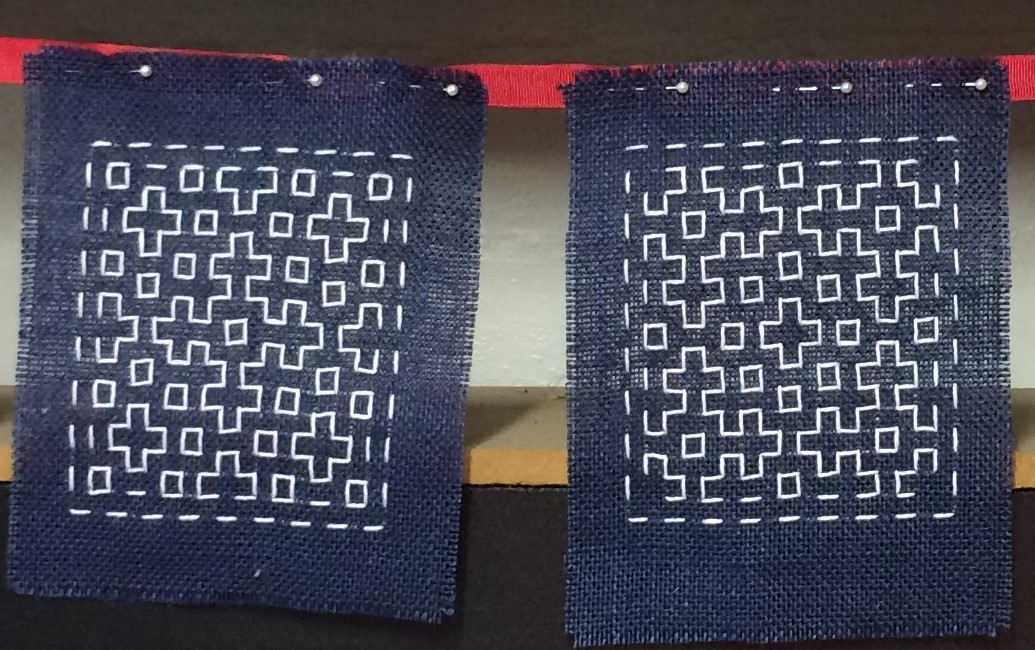}}}
\caption{\label{order3} Two flags featuring Fibonacci snowflakes, from a strip of fractal bunting. The flag on the right is the dual of the one on the left. }
\end{center}
\end{figure}

A series of stitching experiments led us to a way to create an overall design, which moreover is self-dual. 
One such experiment was conducted while devising an entry for a light-hearted competition to design \textit{fractal bunting}, organised in that first dreary year of the pandemic by The Aperiodical (Steckles, 2020). Part of the bunting created is shown in Figure \ref{order3}. The flag on the left shows at its centre a single Fibonacci snowflake of order 3. The flag on the right is its dual, and parts of four order 3 Fibonacci snowflakes flank a configuration of order 1 and 2 Fibonacci snowflakes, familiar from Figure \ref{labbe}. Looking again at the flag on the left, we see parts of this configuration surrounding the order 3 snowflake; this hints at how the design may be extended to repeat in both directions.

These flags are in fact next in a series of designs of which kuchizashi (the top left hand part of Figure \ref{kuchi}) and the crosses-and-mouths of Figure \ref{juji}(b) are respectively the first and second. The first design (kuchizashi) is comprised only of order one Fibonacci snowflakes (mouths). An encoding of {kuchizashi} is $v=w=11$, where $|v|=2=2P_1$. 

The second design contains both order one and order two Fibonacci snowflakes. The design in Figure \ref{juji}(b) can be encoded using $v=w=1001$, where $|v|=4=2P_2$.

Finally, the third design contains the first three Fibonacci snowflakes of Figure \ref{snow}. The design in Figure \ref{order3} is encoded using $v=w=1000110001$, where $|v|=10=2P_3$. 

We are now ready to propose a sequence of hitomezashi designs that feature the Fibonacci snowflakes. However, given that they are arranged differently on the plane than the tessellations of Blondin-Mass\'{e} et al. (2010), and because of their four-fold symmetry, and in tribute to hitomezashi, we call these designs \textit{Pell persimmon polyomino patterns}, or just the Pell persimmons.

\subsection{Pell word and Pell persimmons}
  
Recall the recursion relation (\ref{pell}) for the Pell numbers.  In direct analogy to the Fibonacci-type words of equation (\ref{fibword}), we define the \textit{Pell words} for $n \in \mathbb{N}$ by:
 \begin{equation}
    u_n=\overline{u_{n-1}}
 \ \overline{\widetilde{u_{n-2}} }\  u_{n-1}  \label{words}
 \end{equation}
 with $u_0$ being the empty word and  $u_1=1$. Then the word length is $|u_n|=P_n$. We have not located these words elsewhere in the literature.

The first few Pell words are
 \[
u_1=1\quad u_2=0\ 1\quad u_3=10\ 0\ 01\quad u_4=01110\ 01\ 10001
\]
(Note: small spaces have been used only to allow the eye to see the construction of the word.) The Pell words with odd index are palindromes (that is, $\widetilde{u_{2k+1}}=u_{2k+1}$) and the words with even index are antipalindromes (that is, $\widetilde{u_{2k}}=\overline{u_{2k}}$), which can be proved by induction.

 Let the \textit{Pell persimmon polyomino pattern of order $n$} be the hitomezashi pattern generated by 
 \[
v=w= u_n \widetilde{u_n}.
 \]
 This construction has the `clean' feature that, unlike the formation of the Fibonacci snowflakes from a subset of words given by the somewhat clumsy hybrid expression (\ref{fibword}), there is one definition for all $n$ and all words are used. The recursive definition (\ref{words}) causes the  loops of lower order to appear within and between instances of the largest one in the overall design. Self-duality follows from the following observations:

For $n=2k+1$, since $u_{2k+1}$ is a palindrome,
\[
v=u_{2k+1}\widetilde{u_{2k+1}}=u_{2k+1}u_{2k+1}.
\]
The dual corresponds to swapping zeroes and ones. On the reverse, the correct coincidence of rows and columns to replicate the front fabric recurs beginning at a position shifted vertically and horizontally by $|u_{2k+1}|=P_{2k+1}$ stitches (a number shown previously to be congruent to 1 modulo 4). Recall from Section \ref{code} that reading the encoding an odd number of stitches from the edge effects the $0\leftrightarrow 1$ exchange.

For $n=2k$, since $u_{2k}$ is an antipalindrome,
\[
v=u_{2k}\widetilde{u_{2k}}=u_{2k}\overline{u_{2k}}.
\] 
The pattern on the reverse is rendered by concatenating $\overline{v}$
\[
\overline{v}\, \overline{v} \ldots \overline{v}
=\overline{u_{2k}}\, u_{2k}\,  \overline{u_{2k}}\, u_{2k}
\ldots \overline{u_{2k}}\, u_{2k}=\overline{u_{2k}}\, v\, v \ldots u_{2k}
\]
giving self-duality, up to a shift by $|u_{2k}|$, again both vertically and horizontally.

 We have not found these designs apart from those based only on $u_1$ and $u_2$ in stitch dictionaries. While it is gratifying that these appear new, it would be equally fascinating if they were to be found on garments or homewares in museums. 
 
\begin{figure}
\begin{center}
{{\includegraphics[height=15cm]{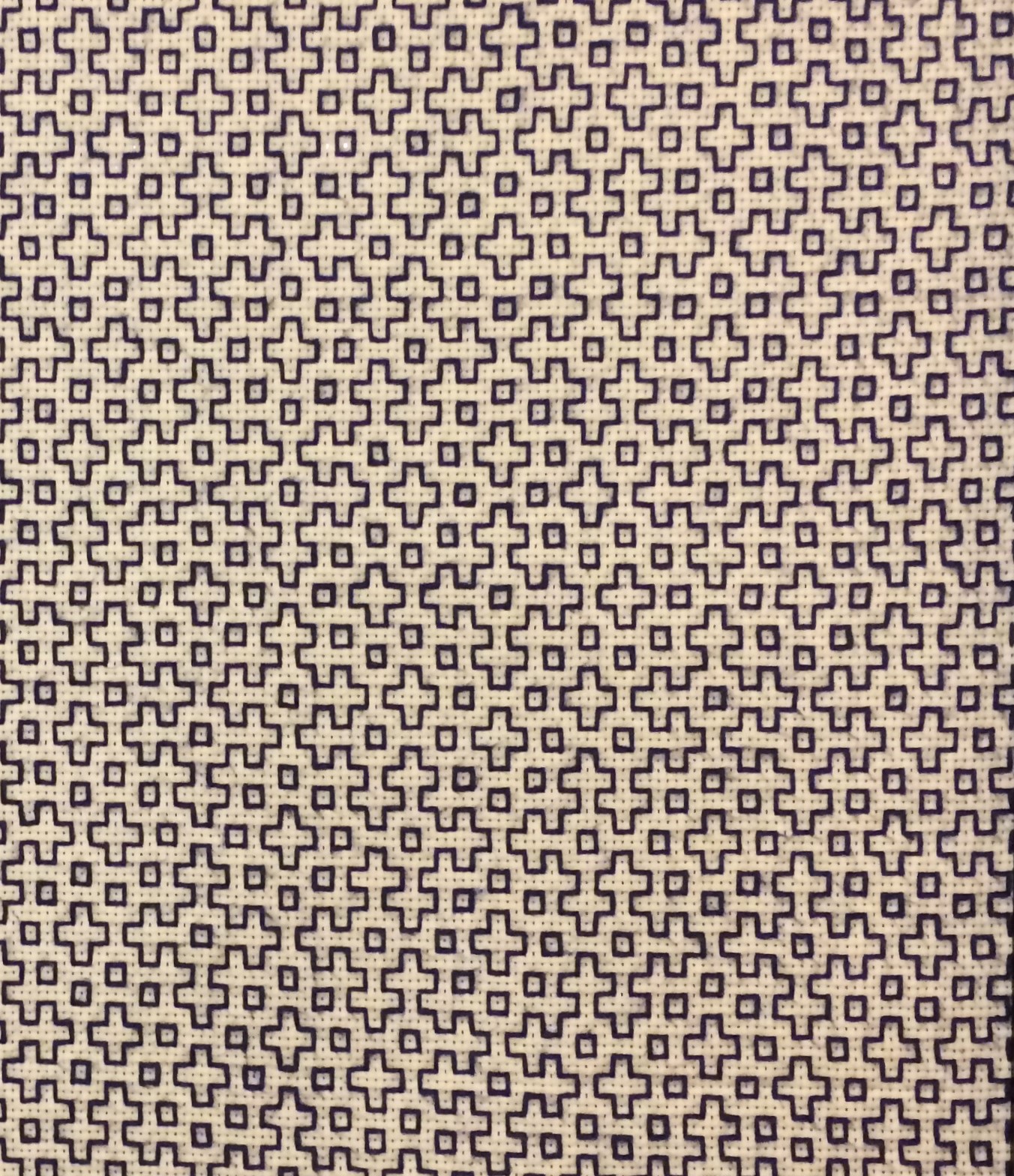}}}
\caption{\label{order5} A single Fibonacci snowflake of order 5 extends from the top to bottom and left to right of this photo. Like Figure \ref{part}, it is a detail from a piece (60 by 84 cm) \textit{Cross and Crown} stitched in 2022. }
\end{center}
\end{figure}

The design in the top left of Figure \ref{kuchi}, Figure \ref{juji}(b), and Figure \ref{order3} provide examples of this construction for order 1, 2 and 3.   Figure \ref{part} shows the Fibonacci snowflake of order 4 and its nested and flanking lower order snowflakes. In fact, it is part of a large piece generated using $u_5$; a different part of this piece, which cannot be shown in its entirety within the confine of a journal page, is shown in Figure \ref{order5}. The loops in the order 5 design were constructed solely using equation (\ref{words}) and only compared after completion to images of the $n=4,\, 5$ Fibonacci snowflakes in the paper of  Blondin-Mass\'{e} et al. (2009). However, the astute reader will realise that we have not given a proof connecting the two constructions in generality. We are confident, however, to make the following conjecture:
 
 \textbf{Persimmon-Snowflake Conjecture:} The largest polyomino in the Pell persimmon polyomino pattern of order $n$ is the Fibonacci snowflake of order $n$. 

\section{Conclusion}
 In future work, the proof of the Persimmon-Snowflake Conjecture, mostly likely adapting methods of Defant and Kravitz (2022), and exploiting  relationships between Pell and Fibonacci numbers, and duality, should be completed. This will form an entire, and technical, paper in its own right. 

\begin{figure}
\begin{center}
{{\includegraphics[height=3cm]{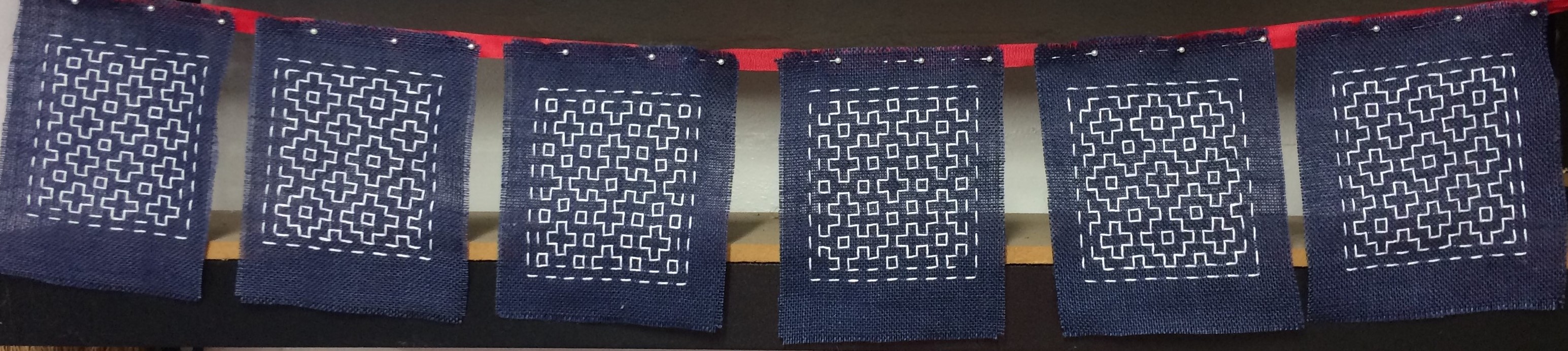}}}
\caption{\label{bunting} The whole piece of fractal bunting. The left-most flag and the two in the centre feature Pell persimmon polynomial patterns as defined in this paper. The remaining three feature generalisations to be explored further in future. }
\end{center}
\end{figure}

Another project, which is currently underway, is to encode hitomezashi designs related to the generalised Fibonacci snowflakes of Blondin-Mass\'{e} et al. (2010) and Ram\'{i}rez et al. (2014), using generalisations of the Pell numbers and hence also of the Pell words. Some of these generalised designs can be seen on the other bunting flags in Figure \ref{bunting}. Some were also used in objects illustrating a previous paper (Seaton, 2021).

Patterns on the triangular lattice (or equivalently, the square lattice with diagonal stitches in one orientation permitted) are described by Defant et al. (2022) as `quite mysterious and difficult to analyse'. The inherent mathematics of the hitomezahi patterns which include diagonal stitches, as well as the other constrained, traditional forms of sashiko described in the introduction, could be considered.

There is also `future work' for you our readers. This article has deliberately been illustrated with examples showing ideas for items to stitch with hitomezashi designs, in a variety of materials.  While drawing is an option, we encourage you to actually stitch them so that you will experience the mathematical concept of duality. Drawing does still promote mindfulness; we lament that those who rush to code the patterns do not experience this benefit of sashiko (Hayes, 2019; Iiduka 2019/2020).

Hitomezashi, although of long heritage, is something of a late-comer in receiving attention from the mathematical fibre arts perspective, but we trust we have convinced you that it was worth the wait.

\section*{Acknowledgements} The image in Figure \ref{labbe} is reproduced with the kind permission of S\'{e}bastien Labb\'{e} from his blog (2010). The images in Figures 4 and 17 have appeared previously in the workshop paper by Hayes and Seaton (2020), those in Figures 7(a), 8, 9(b), 11(a), (b) appeared in an article by Seaton (2021), and the photo in Figure 21 (of which Figure 19 is a detail) appeared at Steckles (2020). All other stitching was created and photographed by KAS. She also photographed the persimmon. Thanks to Chris Taylor and Emerald King, as well as the anonymous reviewers and the Special Issue editors, for their interest and helpful comments.

\section*{Disclosure Statement} The authors have no conflicts to disclose.

\end{document}